\documentclass[a4paper]{article}
\usepackage[english]{babel}
\usepackage[utf8]{inputenc}
\usepackage{amsmath,amsthm}
\usepackage{amsfonts,amssymb}
\usepackage{hyperref}

\usepackage{tikz,tikz-cd} 

\usepackage{xifthen} 

\newcommand{\cat}{\mathcal{C}}

\newcommand{\toposdefaut}{0}
\newcommand{\topos}[1][\toposdefaut]{ 
\ifthenelse{\equal{#1}{0}}{ \mathcal{T} }
{
\ifthenelse{\equal{#1}{1}}{ \mathcal{E} }{ #1 }
}
}

\newcommand{\un}[1][\toposdefaut]{ {1_{\topos[#1]}} }
\newcommand{\soc}[1][\toposdefaut]{ { \Omega_{\topos[#1]} } }
\newcommand{\powerob}[1][]{\mathcal{P}_{\topos[#1]}}
\newcommand{\NN}[1][\toposdefaut]{\mathbb{N}_{\topos[#1] }}
\newcommand{\Zt}[1][\toposdefaut]{\mathbb{Z}_{\topos[#1] }}
\newcommand{\Qt}[1][\toposdefaut]{\mathbb{Q}_{\topos[#1] }}

\newcommand{\Rlscp}[1][\toposdefaut]{\mathbb{R}^{lsc +}_{\topos[#1] }}

\newcommand{\Rt}[1][\toposdefaut]{\mathbb{R}_{\topos[#1] }}
\newcommand{\Ct}[1][\toposdefaut]{\mathbb{C}_{\topos[#1] }}

\newcommand{\Cred}{\mathcal{C}^*_\text{red}}
\newcommand{\Cmax}{\mathcal{C}^*_\text{max}}

\newcommand{\sh}{\textsf{Sh}}

\newcommand{\sub}{\textsf{Sub}}
\newcommand{\set}{\text{Set}}

\newcommand{\slat}{\textsf{sl}}
\newcommand{\module}{\textsf{Mod}}
\newcommand{\rel}{\textsf{Rel}}
\newcommand{\proj}{\textsf{Proj}}
\newcommand{\bil}{\textsf{Bil}}
\newcommand{\site}{\text{Site}}

\newcommand{\hgacoef}[3]{ \begin{pmatrix} #2,#3 \\ #1 \end{pmatrix} }

\newcommand{\scal}[2]{ \left\langle #1 , #2 \right\rangle }

\newcommand{\Qrel}[3][\approx]{ \left[  {#2} #1 {#3} \right] }
\newcommand{\Qfun}[3]{\left[ #1 \approx #2(#3) \right] }

\usepackage{titlesec}
\titleformat{\subsubsection}[runin]{\normalfont}{\thesubsubsection}{0pt}{}[.]

\newcommand{\block}[1]
{
\par \subsubsection{} #1

\bigskip}

\newcommand{\blockn}[1]{\par #1 \bigskip}

\newcommand{\Th}[1]
	{
	\bigskip	
	\textbf{Theorem : }{\itshape #1}
		
	\bigskip
	}

\newcommand{\Prop}[1]
	{

	\bigskip
	
	\textbf{Proposition : }{\itshape #1}
		
	\bigskip
	
	}

\newcommand{\Cor}[1]
	{

	\bigskip
	
	\textbf{Corollary : }{\itshape #1}	
		
	\bigskip

	}

\newcommand{\Lem}[1]
	{

	\bigskip
	
	\textbf{Lemma : }{\itshape #1}
		
	\bigskip
	
	}

\newcommand{\Def}[1]
	{
	
	\bigskip
	
	\textbf{Definition : }{\itshape #1}
	
	\bigskip
	
	}

\newcommand{\Dem}[1]{
	
	\smallskip
	
	\textbf{Proof : } \par
	 {#1} $\square$
	 
	 \bigskip
}

\begin{document}

\pagestyle{plain}
\title{Toposes, quantales and $C^*$ algebras \\ in the atomic case}
\author{Simon Henry}

\maketitle

\renewcommand{\thefootnote}{\fnsymbol{footnote}} 
\renewcommand{\thefootnote}{\arabic{footnote}}

\begin{abstract}
We start by reviewing the relation between toposes and Grothendieck quantales. We improve results of previous work on this relation by giving both a characterisation of the map from the tensor product of two internal sup-lattices to another sup-lattice and a description of the category of internal locales of a topos in terms of the associated Grothendieck quantale. We then construct a convolution product, corresponding to internal composition of matrices, on the set of positive lower semi-continuous functions on the underlying locale of the quantale attached to a topos. In good cases, this convolution product does restrict into a well defined convolution product on a subset of the set of continuous functions and defines a convolution $C^*$ algebra attached to the quantale.
In the last part of this article we investigate in details these attached $C^*$ algebras in the special case of an atomic topos. In this situation the related Grothendieck quantale corresponds to a hypergroupoid. Relatively simple finiteness conditions on this hypergroupoid appear in order to obtain an interesting $C^*$ algebra. This algebra corresponds to a hypergroupoid algebra which comes endowed with an arithmetic sub-algebra and a time evolution. We conclude by showing that the existence of a hypergroupoid satisfying all the requirements attached to a specified atomic topos is equivalent to the fact that the topos is locally decidable and locally separated. Also in this situation the time evolution only depends on the topos and is described by a (canonical) principal $\mathbb{Q}_+^*$ bundle on the topos. The BC-system and more generally the double cosets algebras are special cases of this situation.
\end{abstract}

\tableofcontents

\section{Introduction}

\blockn{Both $C^*$ algebras and toposes yield natural generalizations of topological spaces: locally compact topological spaces correspond to commutative $C^*$ algebras of continuous functions, and sober topological spaces correspond to toposes of sheaves over them. Non commutative $C^*$ algebras and general toposes thus extend the notion of topological space beyond its classical framework. It is hence natural to ask whether these two generalizations are in any sense related to each other, and such a relationship should be extremely interesting for both areas: toposes are  closely related to geometry thanks to their relation with localic groupoids developed in \cite{joyal1984extension}, \cite{bunge1990descent}, \cite{moerdijk1990homotopy} , while non commutative $C^*$ algebras are also intimately related to groupoids (\cite{connes1994book})
and  have proven themselves extremely powerful through their connection with the quantum formalism, the theory of time evolution on type III factors and quantum statistical mechanics.}

\blockn{More recently, a third family of objects appeared in this picture: Quantale. In operator algebra they have been introduced by C.J.Mulvey (see \cite{mulvey1986suppl} and \cite{j2001quantisation}) in an attempt to formalize the notion of ``quantum topology'' studied by R.Giles and H.Kummer in \cite{giles1971topology} and C.A.Akemann in \cite{akemann1969weierstrass}. In topos theory they arise in the description of the category of sup-lattices of a given topos studied in \cite{joyal1984extension} and,  because of the result  \cite{heymans2012grothendieck} (see also the first part of the present article) they completely describe a topos in the sense that a topos endowed with a bound is essentially the same thing as a special kind of quantale, called a ``Grothendieck quantale''.

The type of quantale appearing in operator algebra and in topos theory have extremely different properties: the Grothendieck quantales are quantales of relations on a bound of a topos, and behave like the quantale of relations on a set, in particular they are distributive and modular. On the other side, quantales appearing in operator algebra have different properties, they are in general not modular and correspond to particular subquantales of the  quantale of projections in a Hilbert space, hence deserving the name ``quantum'' in a more precise manner. These differences preclude a straightforward comparison of the  theories of Grothendieck toposes and of operator algebras through the associated  quantales and show that the relation between the two theories is necessarily more involved. We nevertheless use the Grothendieck quantale associated to a topos as a starting point and show in this paper that under suitable hypothesis a Grothendieck quantale can be used to construct a convolution $C^*$ algebra attached to a topos.}

\blockn{We would also like to stress out that modular quantales and Grothendieck quantales are extremely good candidates to be thought of as characteristic one operator algebras. First there are several formal similitudes: the fact that sup-lattices enriched categories are a form of ``characteristic one additive categories", the presence and the important role of the $*$ involution, and other more specific points like the fact that the initial and terminal support of $a$ are given by $a^*a$ and $aa^*$. Secondly, Grothendieck quantales (and conjecturally modular quantales) are interpreted as quantales of relations on objects of topos (see \ref{QSet_Grothrepresentation} and \ref{ModQ_conj}), i.e. as characteristic one matrix algebras. Hence results of sections $3$ underline the fact that there is a close relation between topos theory and  non-commutative geometry in characteristic one. It might be interesting to make this relation more precise, for example by giving an interpretation of the distributivity (Q3) and the modularity (Q9) conditions in term of characteristic one semirings.
}

\blockn{One of the most powerful features of topos theory is probably the internal logic: it has been shown in \cite{boileau1981joyal} that a topos can be considered as a ``mathematical universe'' which differs from the classical universe of sets we are used to work with by the fact that the law of excluded middle and the axiom of choice might fail. This means that there is a way to interpret any mathematical (formally written) statement about sets as a statement about objects of a topos, and any theorem which has a constructive\footnote{By constructive we mean which does not use the law of excluded middle or the axiom of choice.} proof becomes a theorem about objects in any topos. One can consult the part $D$ of \cite{sketches} for a detailed presentation of the logical aspect of topos theory.}

\blockn{We will make an extensive use of this internal logic, and we hope that all the transition between the world of usual sets and the world of sheaves will be clear.
Also, all the mathematics presented in this article are constructive. This is not because of any form of belief from the author that the law of excluded middle should be systematically avoided, but because it was possible to dispense from it without adding too much complexity and it opens the possibility of applying our results internally.}

\blockn{In section $3$ we focus on the relation between ``Grothendieck quantale'' and Grothendieck toposes. Most results of the section $3.1$ to $3.5$ are already well known and present in \cite{pitts1988applications} or \cite{heymans2012grothendieck}. The only originality of our approach is that we give a direct proof that the category of sup-lattices of a topos is the category of modules over a quantale of relations, and then we use this to describe the objects of the topos in terms of this quantale, the previous approach (mainly \cite{pitts1988applications}) generally works in the other direction. In \ref{quantaleintro} we review the basic theory of sup-lattices internal to a topos. In \ref{Slcat} the basic theory of sup-lattices enriched category and we give a characterisation of the sup-lattice enriched categories of modules over a unital quantale. In \ref{Qset} we explain the correspondence between Grothendieck toposes and Grothendieck quantales, and give a description of a topos attached to the quantale $Q$ in term of a notion of $Q$-set completely similar to the notion of $\mathcal{L}$-set when $\mathcal{L}$ is a frame. Section \ref{RelRep} produces a description of the topos attached to a quantale as a classifying topos making the previous correspondence functorial.

The main (and essentially only) new contribution of section $3$ is in \ref{SlatBil} and consists of an extremely elegant description of the locales internal to a topos in terms of ``modular actions'' of the corresponding quantale on classical locales, as well as more generally a description of bi-linear maps between sup-lattices in terms of the corresponding $Q$-modules.

In \ref{ConvProd} we explain why we think that attaching a Grothendieck quantale to a topos is an interesting step towards the construction of $C^*$-algebras.}

\blockn{In section $4$ we focus on the case of an atomic topos; we show that in this case the attached Grothencieck quantale corresponds to a ``Hypergroupoid'' and that under some reasonable finiteness assumptions there is indeed a ``Hypergroupoid $C^*$-algebra'' attached to that quantale in the way sketeched in \ref{ConvProd}. This $C^*$ algebra comes in two forms: a reduced algebra and a maximal algebra; in both cases it comes with a natural and explicit time evolution attached through Tomita theory to a ``regular'' representation, and with a generating $\mathbb{Z}$ sub-algebra with interesting combinatorial properties. We also characterise in section \ref{septop} the atomic toposes for which the construction is possible as the locally decidable locally separated toposes. We also show that in this situation the time evolution is canonical and described by a principal $\mathbb{Q}_+^*$ bundle. The main example of this situation are the well-known double-cosets algebras.}

\section{Notations and preliminary}

\blockn{By toposes we always mean Grothendieck toposes, i.e. categories of sheaves on a Grothendieck site. Most of the basic notions of topos theory can be found in \cite{maclane1992sheaves}, for the others we will give precise references in \cite{sketches}.
\begin{itemize}

\item if $\cat$ is a category (or a topos) then $|\cat|$ denotes the set (or the class) of objects of $\cat$. The symbol $\cat$ denotes the set (or class) of all maps, and we will equivalently use the notation $\cat(a,b)$, $\hom(a,b)$ or $\hom_{\cat}(a,b)$ for the set of morphisms from $a$ to $b$.

\item The letter $\topos$ always denotes a topos.

\item $\soc$ denotes its sub-object classifier.

\item $\un$ denotes its terminal object.

\item If $X \in |\topos|$ then $\powerob(X)$ stand for the power object of $X$ (isomorphic to $\soc^X$), and $\sub(X)$ for the set of sub-objects of $X$, i.e. the set of global sections of $\powerob(X)$.

\item A set (or an object $X \in |\topos|$) is said to be inhabited if it satisfy (internally) $\exists x \in X$ (which in constructive mathematics is stronger than the assertion that $X$ is not empty). For an object of a topos it is equivalent to the fact that the canonical map $X \rightarrow 1$ is an epimorphism.

\item $\NN$, $\Zt$ and $\Qt$ denote the sheaves of natural numbers, integers and rational numbers in the topos $\topos$. As we restrict ourselves to Grothendieck toposes, they are simply $p^*(\mathbb{N})$,$p^*(\mathbb{Z})$ and $p^*(\mathbb{Q})$ where $p$ is the canonical geometric morphism from $\topos$ to the topos of sets.

\item $\Rt$ denotes the sheaf of continuous real numbers, i.e. two sided Dedekind cuts (see \cite{sketches} D4.7). It can be described externally by the following properties: for any $X \in |\topos|$, $\hom(X,\Rt)$ is the set of continuous functions from the underlying locale of $X$( whose frame of opens is $\sub(X)$) to the space of real numbers\footnote{In a non-boolean context, the ``space of or real numbers'' has to be interpreted as ``the formal locale of real numbers''.}.

\item $\Rlscp$ denotes the set of positive lower semi-continuous real numbers (possibly infinite). In presence of the law of excluded middle it is the set $\mathbb{R}_+ \cup \{ \infty\}$. In a topos it is the sheaf defined by the fact that $\hom(X,\Rlscp)$ is the set of functions from the locale $\sub(X)$ to $\mathbb{R}_+ \cup \{ \infty\}$\footnote{Here we have assumed the law of excluded middle in the topos of set in order to simplify the notation.} endowed with the topology where the $(a,\infty]$ are a basis of open sets, i.e. it is the set of lower semi-continuous functions (possibly infinite) on the locale $\sub(X)$. Internally, $\Rlscp$ is defined as the set of $P \subset \Qt$ such that if $q<0$ then $q \in P$, and $q \in P \Leftrightarrow \exists q' \in P, q<q' $. See \cite{sketches} D4.7.

\item A sub-quotient of an object $X \in |\topos|$ is a quotient of a sub-object of $X$ (or equivalently, but less naturally, a sub-object of a quotient).

\item A proposition (internal to a topos) is said to be decidable if it is complemented (i.e. equal to its double negation). An object is said to have decidable equality, or to be decidable if its diagonal embedding $X \rightarrow X \times X$ is complemented.

\end{itemize}}

\blockn{Also an object $B \in |\topos|$ is said to be a {\it bound} of $\topos$ if any object of $\topos$ can be written as a sub-quotient of an arbitrary co-product of copies of $B$ (see \cite{sketches} B3.1.7.). Equivalently $B$ is a bound of $\topos$ if $\sub(B)$ is a generating family of $\topos$, i.e. $\sub(B)$, seen as a full subcategory of $\topos$ and endowed with the restriction of the canonical topology of $\topos$, forms a site of definition for $\topos$.
This means essentially that $B$ is big enough to generate $\topos$: in the topos of $G-\set$ for a group $G$, an object $X$ is a bound if and only if the map $G \rightarrow Aut(X)$ is injective. A topos is the topos of sheaves over a locale if and only if $\un$ is a bound (see \cite{sketches}Definition A4.6.1 and theorem C1.4.6). When a topos is given by a site, the simplest way to obtain a bound is to choose an object which contains a copy of each representable object (for example, the direct sum of all the representable objects, see \cite{sketches}B3.1.8(b)).

Existence of a bound, together with the existence of enough (co)limits characterize Grothendieck topos among elementary topos (see \cite{sketches} C2.2.8)}

\blockn{A set $X$ (or an object of a topos) will be said to be finite if it is (internally) Kuratowski finite, i.e. if $\exists n \in \mathbb{N}$, $x_1, \dots, x_n \in X$ such that for all $x \in X$ $\exists i, x=x_i$. On can consult \cite{sketches} D4.5 for the theory of Kuratowski finite set.
Roughly, a quotient of a finite set is finite, but the proof that a subset of a finite set is finite requires the subset to be complemented and may fail in full generality. If a set $X$ is finite and has decidable equality, then there exists $n \in \mathbb{N}$ such that $X$ is isomorphic (internally\footnote{As the isomorphism is not canonic, it might not lift to a global isomorphism if we are working internally in a topos.}) to $\{1, \dots, n \}$, and a subset of $X$ is finite if and only if it is complemented.
}

\section{Toposes, quantales and sup-lattices}

\subsection{Introduction}
\label{quantaleintro}
\block{Let $X$ be any object of a topos $\topos$, we will denote by $\rel(X)$ the set of relations on $X$, i.e. the set $\sub(X \times X)$ of sub-objects of $X \times X$. Then $Q=\rel(X)$ is endowed with several structures: 
\begin{enumerate}
\item[(Q1)] The inclusion of subobjects gives an order relation on $Q$.

\item[(Q2)] $Q$ has arbitrary supremum for this order relation.

\item[(Q3)] Finite intersections distributes over arbitrary union: $a \wedge \bigvee_i b_i = \bigvee_i (a \wedge b_i)$.

\item[(Q4)] There is an associative composition law on $Q$ defined internally by $R P = \{(x,y) | \exists z \in X, x R z \text{ and } z P y \}$.

\item[(Q5)] The composition law is order preserving and commutes to supremum in each variable.

\item[(Q6)] The diagonal subobject of $X$ provide an element $1 \in Q$ which is a unit for the composition law.

\item[(Q7)] There is an order preserving involution: $R \mapsto R^* = \{ (y,x) | x R y \} $ of $Q$.

\item[(Q8)] For all $x,y \in Q$ one has: $(x y) ^* = y^* x ^*$.

\item[(Q9)] For all $x,y,z \in Q$ one has $x \wedge y z \leqslant y ( y^* x \wedge z)$.

\end{enumerate}

If we assume additionally that $X$ is a bound of $\topos$, then one has additionally:
 \begin{enumerate}
 
\item[(Q10)] For all $q \in Q$ there are two families $(v_i)_{i \in I},(u_i)_{i \in I}$ of elements of $Q$ such that: $\forall i, u_i u_i^* \leqslant 1, v_i v_i^* \leqslant 1$ and $\top = \bigvee_i v_i u_i^*$. where $\top$ denote the top element of $Q$.

\end{enumerate}

Some of these points deserve a proof and a few comments.

\begin{itemize}

\item $(Q9)$ is called the modularity law. 
It is easy to prove using internal logic:
Let $(a,b) \in (x \wedge y z)$. One has: $(a,b) \in x$ and there exists $c \in X$ such that $(a,c) \in y$ and $(c,b) \in z$. Hence $(a,c) \in y$ and $(c,b) \in (y^* x \wedge z)$ so $(a,b) \in y(y^*x \wedge z)$. As this proof uses only intuitionist logic, it is valid in any topos.

\item $(Q3)$ is sometimes also called the modularity law, which gives rise to a conflict of terminologies. We will prefer the term distributivity law for $(Q3)$.

\item $(Q10)$ expresses the fact that, as $X$ is a bound of $\topos$, $X \times X$ has to be a sub-quotient of a co-product of a set $I$ of copies of $X$. 

Indeed, in this situation, there is a family $(u_i,v_i)_{i \in I}$ of partial applications from $X$ to $X \times X$. A partial application $f$ from $X$ to $X$ can be represented by its graph: the relation $R$ such that $(y R x)$ if and only if $f(x)$ is defined and $y = f(x)$. A relation $R$ on $X$ is the graph of a partial application if and only if $R R^* \leqslant 1$.
So one has two families of relations on $X$, also denoted $(u_i)$ and $(v_i)$, such that for all $i$, $u_i u_i^* \leqslant 1$ and $ v_i v_i^* \leqslant 1$. The relation $\bigvee_i v_i u_i^*$ is the union of the image of $X$ in $ X \times X$ by all the partial maps $(v_i,u_i)$. So the relation $\bigvee_i v_i u_i^* = \top$ expresses the fact that the corresponding map is onto.

\end{itemize}

}

\block{\label{quantaleintro_def}\Def{A Set satisfying $(Q1)$ and $(Q2)$ is a sup-lattice. A Set satisfying $(Q1)$, $(Q2)$, $(Q4)$ and $(Q5)$ is called a quantale, (unital if it also satisfies $(Q6)$).
We will call a modular quantale, a quantale satisfying all the axioms from $(Q1)$ to $(Q9)$, and a Grothendieck quantale one satisfying all the axioms from $(Q1)$ to $(Q10)$.
}

The term quantale is due to C.J. Mulvey in \cite{mulvey1986suppl}. The name Grothendieck quantale has been given by H.Heymans in \cite{heymans2010sheaves}. For the term ``modular quantale'', our terminology differs slightly from previous work (like  \cite{heymans2012grothendieck}), where the axiom $(Q3)$ is not included in the definition of a modular quantale. The main reasons for our choice of terminology is simply that we only want to consider quantale that arise as relations on objects in a topos and hence satisfy the axiom $(Q3)$ and also that we think it is more natural to assume a compatibility between intersection and supremum (given by $(Q3)$) as soon as we assume both a compatibility between intersection and the composition law (given by $(Q9)$) and a compatibility between the composition law and supremum (given by $(Q5)$).

}

\block{The main result relating toposes to quantale (which should probably be attributed to Frey and Scedrov in \cite{Freyd_Alegories}), is the fact that if $\topos$ is a topos and $B$ is a bound of $\topos$ then $\topos$ can be completely reconstructed from the Grothendieck quantale $Q = \rel(B)$. And that every Grothendieck quantale can be written (essentially uniquely) in the form $\rel(B)$ for a bound $B$ of a topos $\topos$.

This result (at least its first part, the second part being a little harder) can actually be proven directly using the following construction:

\Def{If $Q$ is a Grothendieck quantale, we will denote by $\site(Q)$ the site whose objects are the $q \in Q$ such that $q \leqslant 1$, whose morphisms are given by:
 \[ \hom(q,q') = \{f \in Q | 1 \wedge f^* f = q \text{ and } f f^* \leqslant q' \}. \]
The composition is given by $f \circ g =fg$. The identity morphism of $q$ is $q$ itself. And a Sieve $J$ on an object $q$ is covering if:

\[ \bigvee_{q' \leqslant 1, f \in J(q')} f f^* = q \] 

}

The fact that for any Grothendieck quantale $\site(Q)$ is indeed a site is not straightforward. Apparently\footnote{We checked it, but unfortunately, it does not seems that a proof of this kind had ever been published.} it can be checked directly, but this proof is quite long and is not necessary because one has a more abstract proof, using the following easier proposition, and theorem \ref{QSet_Grothrepresentation}.

\Prop{If $Q = \rel(B)$, for a bound $B$ of a topos $\topos$, then $\site(Q)$ is the site of subobjects of $B$. In particular, it is a site of definition for $\topos$. }

\Dem{We use the same kind of argument as the proof that $Q$ satisfies $(Q10)$.

As $1\in Q = \sub(B \times B)$ corresponds to the diagonal sub-object of $B \times B$, an element $q \in Q$ such that $q\leqslant 1$ corresponds to a unique subobject of $B$.
Let $q$ and $q'$ be two subobjects of $B$, and $f \in Q =\rel(B)$ satisfying the two conditions $1\wedge f^*f=q$ and $ff^*\leqslant q'$. The first condition asserts (internally) that all $x$ such that $\exists y$ , $y f x$ is an element of $q$, and the second condition asserts that if $y fx$ and $y' f x$ then $y=y'$ and $y\in q'$. This is exactly the condition that characterises the graph of a function from $q$ to $q'$, hence $\hom_{\site(Q)}(q,q')$ is indeed isomorphic to $\hom_{\topos}(q,q')$ and as the composition of relations extends the composition of functions this correspondence is indeed an equivalence of categories.

It only remains to check that the topology of $\site(Q)$ is indeed the canonical topology of the topos, but for any collection of map $f_i :q_i \rightarrow q$, the sub-object $f_i f_i^* \leqslant q$ is exactly the image of $f_i$ in $q$ and hence the condition that:

\[\bigvee_{i} f_i f_i^* = q\]
simply asserts that the family is jointly surjective.
}

One of our goal is to provide a way to reconstruct $\topos$ from $Q$ without using sites. This construction gives an alternative to sites for working with Grothendieck topos.

}

\subsection{The category $\slat(\topos)$ of sup-lattices}
\label{slat}

\blockn{In this subsection we recall the definition and basic properties of the categories of sup-lattices of a topos as it is defined and studied in \cite{joyal1984extension}. We will not give any proofs, but most of them are straightforward and they all can be found in \cite{joyal1984extension}.}

\block{\Def{We will denote by $\slat(\topos)$ the category of sup-lattices of $\topos$. This means that objects of $\slat(\topos)$ are objects $S$ of $\topos$ endowed with a relation $(\leqslant) \subset S \times S$, such that $(S,\leqslant)$ is internally a sup-lattice (i.e. $\leqslant$ is a partial order relation and $S$ admits arbitrary supremum). The morphisms are the morphisms $f:S \rightarrow S'$ in $\topos$ which (internally) preserve supremums (and hence also the order relation).}

In all of this sub-section, we will fix the topos $\topos$ and work internally inside it. We will denote simply by $\slat$ the category $\slat(\topos)$ and consider objects of $\topos$ as usual sets.

}

\block{\label{slat_involution}Although we use the term ``sup''-latices, it is a classical fact of ordered sets theory that if every subset admits a supremum then every subset also admits an infimum, and hence a sup-lattice is the same thing as an inf-lattice. The term ``sup'' is here to denote the fact that we are considering sup-preserving morphisms (which are different from inf-preserving morphism).

This duality has a consequence on the category $\slat$: it is endowed with an involutive contravariant functor, that we will denote by $(\_)^*$. Indeed if $S$ is a sup-lattice then if we define $S^*$ as being $S$ endowed with the reverse order relations it is again a sup-lattice, and if $f$ is a morphism then we denote by $f^*$ its right adjoint (it always exists because $f$ commutes to supremum) and as $f^*$ commutes to infimum it is a morphism of sup-lattice for the opposite order relations. One has $f^{**} = f$ because of the reversing of the order relations, and hence $*$ is an involutive anti-equivalence of categories. 

\bigskip

This involution allows to compute colimits in the category of sup-lattices: indeed one can easily check that $\slat$ has all limits and that they are computed at the level of the underlying set. As $(\_)^*$ transforms co-limits into limits, the computation of a colimit can be brought to the computation of a limit.

}

\block{If $X$ is a set, then $\powerob(X)= \soc^X$ (the power object of $X$) is a free sup-lattice generated by $X$, i.e.:

\[ \hom_{\topos}(X,S) = \hom_{\slat}(\powerob(X),S) \]

This adjunction formula turns $\powerob$ into a functor from sets to $\slat$ that sends a map $f: X \rightarrow Y$ to the direct image map $\powerob(f): \powerob(X) \rightarrow \powerob(Y)$.
}

\block{Knowing how to construct a free sup-lattice (using $\mathcal{P}$) and a quotient of sup-lattice (using the involution $^*$), one can construct sup-lattices by ``generators and relations''. More precisely, if $I$ is a set, and $R$ is a family of couples of subsets $(r_1,r_2)$ of $I$, interpreted as relation of the form: 

\[ \bigvee_{x \in r_1} x \leqslant \bigvee_{y \in r_2} y \]

then the sup-lattice presented by the set of generators $I$ and the set of relations $R$ identifies with:

\[ \{V \subset I | \forall (r_1,r_2) \in R, (r_2 \subset V) \Rightarrow (r_1 \subset V) \} \]

}

\block{If $S$ and $S'$ are sup-lattices, then the set of sup-lattice morphisms between $S$, $S'$ is again a sup-lattice, for the point-wise ordering, with supremum computed point-wise. This sup-lattice is denoted by $[S,S']$.

\bigskip 

These internal $\hom$ objects come with a monoidal structure given by the universal property:

\[ \hom(M\otimes N, P ) \simeq \hom( M , [N,P] ) \] 

Equivalently, the morphisms from $M \otimes N$ to $P$, are the applications from $M \times N$ to $P$ which are morphisms of sup-lattice in each variable (when fixing the other variable). We will call such applications bilinear maps from $M \times N$ to $P$.
The explicit construction of the tensor product is conducted exactly as for modules over a ring by a construction by generator (the $(m \otimes n)$ for $m,n \in M \times N$) and relations expressing the notion of bi-linear map.
}

\block{In addition of being a closed monoidal category endowed with an involution, the category $\slat$ also satisfies the following interesting properties. 

\[ \soc \otimes N = N \]
\[ [\soc,N] = N \]
\[ M^* = [M,\soc^*] \]
\[ (M \otimes N)^* = [M,N^*] \]

In particular, even if we will not use this concept here, this means that $\slat$ (endowed with all these structures) is a $*$-autonomous category in the sense of \cite{barr1979autonomous}, with $\soc^*$ as dualizing object.

 }

\block{\label{slat_functoriality} Let $\topos$ and $\topos[1]$ be two toposes, and $f$ a geometric morphism from $\topos$ to $\topos[1]$. Let also $S$ be a sup-lattice in $\topos$, then $f_*(S)$ is a sup-lattice in $\topos[1]$: indeed (working internally in $\topos[1]$) if $P$ is a subset of $f_*(S)$ then by adjunction there is a map from $f^*(P)$ to $S$, we can consider the supremum $s$ of the image of this map. As $s$ is a uniquely defined element, it is a global section of $S$, i.e. an element of $f_*(S)$. From here one can check that $s$ is also a supremum for $P$. This defines a functor $f_*: \slat(\topos) \rightarrow \slat(\topos[1])$. We also note that $f_*$ preserves bi-linear maps between sup-lattices.

\bigskip

In the other direction, if $S$ is a sup-lattice in $\topos[1]$ then $f^*(S)$ is in general just a pre-order set in $\topos$, but one can construct a completion, denoted by $f^\#(S)$. In order to do so, we chose any presentation by generators and relations of $S$ (for example, taking all elements and all relations), and then we define $f^\#(S)$ by generators and relations using the pull-back of the system of generators and relations chosen for $S$. At first sight, it is not clear that this definition of $f^\#(S)$ does not depend of the presentation of $S$, but one can prove an adjunction formula:

\[ \hom_{\slat(\topos)}(f^\#(S),T) \simeq \hom_{\slat(\topos[1])}(S,f_*(T)) \]

which is natural in $T$. This implies that $f^\#(S)$ does not depend on the presentation of $S$, that it is functorial in $S$ and that $f^\#$ is a left adjoint of $f_*$. We will use the same technique in the proof of the third point $2$ of proposition \ref{Slcat_embedings}. This result can actually be seen as a special case of the first two points of proposition \ref{Slcat_embedings} applied internally in $\topos[1]$ to the category $\cat =\slat(\topos)$ with $B=\soc$.
}

\subsection{Categories enriched over $\slat$}
\label{Slcat}
\block{Thanks to the monoidal structure on $\slat$ one can talk about $\slat$-enriched category. Precisely, a $\slat$-enriched category $\cat$ is a category such that morphism sets are endowed with an order relation which turns them into sup-lattices and composition into a bi-linear map. 

\bigskip

Here are the two main examples of $\slat$-enriched category we want to consider:

\Prop{Let $\topos$ be a Grothendieck topos, then $\slat(\topos)$ is a $\slat$ enriched category. }

\Dem{If $S$ and $S'$ are two objects of $\slat(\topos)$ and $p$ denotes the structural geometric morphism from $\topos$ to the topos of sets, then 

\[ \hom(S,S') = p_*([S,S']) \] 

which is a sup-lattice thanks to \ref{slat_functoriality}. The composition is a bi-linear map because it is given (through an application of $p_*$) by an internal bi-linear map:
 \[ [S,S'] \times [S',S''] \rightarrow [S,S'']. \]
}

\bigskip

A (unital) quantale, as defined in \ref{quantaleintro_def}, is exactly a monoid object of $\slat$, i.e. it is a sup-lattice endowed with the structure of a (unital) monoid such that the composition law is bi-linear. A right (or left) module over a unital quantale $Q$, is a sup-lattice $S$ endowed with a right (or left) action of the underlying monoid of $Q$ such that the corresponding map $S \times Q \rightarrow S$ is bi-linear.

The category of right modules over $Q$ (with $Q$-linear morphisms) is denoted by $\module_Q$, this is the other important example of $\slat$-enriched category we will consider.

\bigskip

If one thinks of the supremum of a family of elements as a form of addition, a sup-lattice enriched category is really close to being an additive category (maybe something we would like to call a ``locally complete characteristic one additive category'' as our addition is characterised by the properties that $x+x=x$ ). The following two results are in this spirit.

}

\block{\label{Slcat_biproduct} From the technique of computation of co-limits in $\slat$ explained in \ref{slat_involution} one can see that the co-product of a family of objects in $\slat$ is isomorphic to the product of the same family. This is actually a general well known\footnote{it appears, for example, under a slightly different form in \cite{Freyd_Alegories} 2.214 and 2.223.} result:

 \Prop{Let $\cat$ be a $\slat$-enriched category, let $(A_i)_{i \in I}$ be a family of objects of $\cat$ and $A$ be an object of $\cat$, then the following three conditions are equivalent:
\begin{enumerate}
\item $A$ is the co-product of the family $(A_i)$.

\item $A$ is the product of the family $(A_i)$.

\item There are two families of morphisms $f_i: A_i \rightarrow A$ and $p_i: A \rightarrow A_i$ such that  $\sup_i f_i \circ p_i = Id_A$ and for all $i,j \in I$:
\[ p_j \circ f_i = \sup \{ f: A_i \rightarrow A_j | i=j \text{ and } f = Id_{A_i} \}
\footnote{If we assume the law of excluded middle, or more specifically that the set of indicies $I$ has a decidable equality, then this formula reduces to the more classical: $p_i \circ f_i = Id_{A_i}$ and $p_j \circ f_i = 0$ if $ i \neq j$ } 
\] 

\end{enumerate}

Moreover, in this situation, the morphisms $f_i$ and $p_i$ given in $3.$ are the natural morphisms asserting that $A$ is the (co)-product of the $A_i$.
} 

\Dem{Passing from $\cat$ to $\cat^{op}$ preserves property $3.$ and exchanges properties $1.$ and $2.$, hence it is enough to show that $2.$ and $3.$ are equivalent. 

We will start by showing that $3. \Rightarrow 2.$.

We assume $3.$ holds, in particular $A$ is already endowed with maps ($p_i$) from $A$ to $A_i$ for each $i$, we have to show that $A$ and the $(p_i)$ are universal for this property.

Let $X \in \cat$ be any object and assume we have a collection of map $h_i: X \rightarrow A_i$. 
Let $h = \sup_i (f_i \circ h_i): X \rightarrow A$. Then for every $i$:
\[  p_i \circ h = \sup_j p_i \circ f_j \circ h_j = \sup_j \sup \{ f \circ h_j | i=j \text{ and } f= Id_{M_i} \} = h_i. \]

We also have to show that this map is unique:
let $h'$ be any other map from $X$ to $A$ such that for every $i$, $p_i \circ h' = h_i$. Then:
\[ h= \sup f_i \circ h_i = (\sup f_i \circ p_i) \circ h' = h'.\]

Assume now that $A$ is the product of the $A_i$. The maps $p_i$ are the structural maps, the maps $f_i$ are uniquely defined morphisms (using the universal property of the product) by the formula given for $p_j \circ f_i$. Hence the formula for $p_j \circ f_i$ holds by definition, and the equality $\sup_i f_i \circ p_i = Id_A$ because of the relation

\[ p_j \circ \sup_i f_i \circ p_i = p_j \]

(obtained by the same computation as in the first part of the proof) and the uniqueness in the universal property of the product.

}

This proposition has interesting consequences: First, any $\slat$-enriched functor will automatically preserve each product and each co-product (because $3.$ is clearly preserved by any $\slat$-enriched functor).

Additionally, one can describe the morphisms between two co-products (or products) by something which looks like (infinite) matrix calculus. More precisely, a morphism $f$ from $\coprod_{j \in J} A_j$ to $\coprod_{i \in I} B_i$ is the same thing as a morphism from $\coprod_{j \in J} A_j$ to $\prod_{i \in I} B_i$, hence it is given by the datum of a morphism $f_{i,j}: A_j \rightarrow B_i$ for each $i$ and each $j$.

The composition with a $g: \coprod_{i \in I} B_i \rightarrow \coprod_{k \in K} C_k$, is:

\[ g \circ f = \bigvee_{i \in I} g_{k,i} \circ f{i,j} \]

In the special case where all the $A_j$ and $B_i$ are isomorphic to a same object $A$, then $\hom(A,A) = Q$ is a quantale and we will denote by $M_{I,J}(Q)$ the set of morphisms from $A^{(J)}$ to $A^{(I)}$, which can be identified with $Q^{I \times J}$.

}

\block{\label{Slcat_embedings}The next result can be thought of as a $\slat$-enriched form of the Mitchell embedding theorem which asserts that every abelian category is a full sub-category of a category of modules over a ring, but restricted to the case where there are enough ``projective'' objects.

\Prop{Let $\cat$ be a $\slat$-enriched category, $A$ an object of $\cat$ and $Q = \hom_{\cat}(A,A)$.

\begin{enumerate}

\item $Q$ is a quantale for composition, and $R_A: X \mapsto \hom_{\cat}(A,X)$ induces a functor from $\cat$ to $\module_Q$.

\item If $\cat$ has all co-limits, then $R_A$ has a left adjoint denoted $T_A: Y \mapsto Y \otimes_Q A$.

\item If in addition $R_A$ commutes to co-equaliser (we will say that $A$ is regular projective), then $T_A$ is fully faithful.

\item If in addition $\cat(A,\_)$ detects isomorphisms (i.e. if $f$ is a map in $\cat$ such that $\cat(A,f)$ is an isomorphism then $f$ is an isomorphism), then $R_A$ and $T_A$ realize an equivalence of categories between $\cat$ and $\module_Q$.

\end{enumerate}
 }

\Dem{
\begin{enumerate}
\item Clear: As $\cat$ is an $\slat$ enriched category, composition are bilinear, hence $Q=\hom_{\cat}(A,A)$ is a quantale for composition, the action of $Q$ on $\hom_{\cat}(A,X)$ by pre-composition is also bi-linear, and for any morphism $f:X \rightarrow Y$ the induced morphism from $\hom_{\cat}(A,X)$ to $\hom_{\cat}(A,Y)$ is a $Q$-linear morphism of sup-lattices.

\item Let $X$ be a right $Q$ module, then (in $\module_Q$) one has a surjection $p: \coprod_{x \in X} Q \twoheadrightarrow X$. Let $f_1,g_1: K \rightrightarrows \coprod_{x \in X} Q $ be the kernel pair of $p$. Let $p_2: \coprod_{k \in K} Q \twoheadrightarrow  K $, and let $f = f_1 \circ p_2$ and $g = g_1 \circ p_2$. 

$X$ is the co-equaliser of the two $Q$ linear maps (for the right action) $f$ and $g: \coprod_{k \in K} Q \rightrightarrows  \coprod_{x \in X} Q $, which correspond to elements of $M_{X,K}(Q)$.

Let $A^{(X)} = \coprod_{x \in X} A$ and $A^{(K)} = \coprod_{k \in K} A$. Thanks to a remark done in \ref{Slcat_biproduct}, maps between $A^{(K)}$ and $A^{(X)}$ can also be identified with elements of $M_{X,K}(Q)$, hence there are two maps corresponding to $f$ and $g$ from $A^{(K)}$ to $A^{(X)}$. We define $T_A(X)$ to be the co-equaliser of these two maps.

One easily checks that for any $B \in \cat$, there is a canonical (functorial in $B$) isomorphism $Hom(T_A(X),B) \simeq Hom(X,R_A(B))$ (they are the same once we develop all the inductive limits involved) which implies both the adjunction between $T_A$ and $R_A$ and the functoriality of $T_A$.

\item As $T_A(X)$ is computed as the co-equaliser of two arrows $f,g: A^K \rightrightarrows A^X$ such that the co-equalizer of $R_A(f),R_A(g)$ is $X$, if $R_A$ commutes to co-equalizer then one can deduce that $R_A(T_A(X)) \simeq X$ which (thanks to the adjunction) means that $T_A$ is fully faithful.

\item We already know that $X \simeq R_A\circ T_A(X)$ (by the unit of the adjunction). Let $c_X: T_A(R_A(X)) \rightarrow X$ be the co-unit of the adjunction then, $R_A(c_X): R_A(T_A(R_A(X))) \rightarrow R_A(X)$ is a retraction (by general properties of the unit and co-unit) of the unit of the adjunction evalued in $R_A(X)$ (i.e. the canonical map $R_A(X) \rightarrow R_A(T_A(R_A(X))) $) but this map is known to be an isomorphism, hence $R_A(c_X)$ is an isomorphism and since $R_A$ detects isomorphism, we proved that $c_X$ is an isomorphism.

\end{enumerate}

}

}

The following theorem can then be seen as a corollary of the previous proposition.

\block{\label{Slcat_topos}\Th{Let $\topos$ be a Grothendieck topos, and $B$ a bound of $\topos$. Then $\hom_{\topos}(B, \_)$ induces (one half of) an equivalence of categories from $\slat(\topos)$ to $\module_{Q}$ where $Q$ is the quantale $\rel(B)$.  }

This result is essentially the same as the theorem 5.2 of \cite{pitts1988applications}.

\Dem{We will prove that with $\cat = \slat(\topos)$ and $A= \powerob(B)$, all the hypotheses of the four points of the previous proposition are verified, and $Q = \rel(B)$. 

Note that:
\[ \hom_{\topos}(B, S) = \hom_{\slat(\topos)}(A,S). \]

\begin{itemize}
\item $\slat(\topos)$ has all co-limits (and also all limits) because they can be computed internally in $\topos$.

\item $R_A$ commutes to co-equalizer because of the following formula:

\[ R_A(X) = \hom_{\topos}(B,X) \simeq \hom_{\topos}(B,X^*)^* \simeq \hom(A,X^*)^* \simeq \hom(X,A^*)^* \]

And the last term clearly commutes to every inductive limit.

\item $Q$ is identified with $\rel(B)$ through the isomorphism:

\[ \hom_{\slat(\topos)}(\powerob(B), \powerob(B)) \simeq \hom_{ \topos}(B, \powerob(B)) \simeq \sub(B \times B) \]

Internally, this corresponds to the map which sends a morphism $f$ to the relation $y(R_f)x:= `x  \in f(\{y\})'$. The fact that composition of morphisms coincides with the composition of relations is checked internally:

\[ z R_f R_g x = (\exists y, x \in f(\{y \}) \text{ and } y \in g(\{ z\}) ) = x R_{f \circ g} z. \]

\item $R_A$ detects isomorphisms:

Let $f: S \rightarrow S'$ such that $R_A(f)$ is an isomorphism. For any subobject $U \subset B$, every map $t: U \rightarrow S$ can be extended canonically to a map $\tilde{t}: B \rightarrow S$ by the (internal) formula: 

\[ \tilde{t}(x) = \sup \{ y |  x \in U \text{ and } y =t(x)  \} \]

If $t$ is a map from $B$ to $S$, we can restrict $t$ to $U$ and then extend $t |_U$ into $\tilde{t}$, one has then the formula 
\[ \tilde{t}(x) = \sup \{ y |  x \in U \text{ and } y =t(x)  \} = t. \delta_U\]
where $\delta_U$ is the element of $Q$ corresponding to the diagonal embedding of $U$ in $B \times B$ and the product is the natural right action of $Q$ on $\hom(B,S)$.
Finally, as $\delta_U^2 = \delta_U$, $\hom(U,S)$ is identified with $\hom(B.S).\delta_U$.

As $R_A(f)$ is an isomorphism, all the maps $\hom(U,f)$ for every subobject $U$ of $B$ are isomorphisms, because they are retractions of the map $R_A(f)$. The object $B$ being a bound of $\topos$, the subobjects of $B$ form a generating family and so $f$ is an isomorphism.

\end{itemize}

}

}

\subsection{Quantale Sets}
\label{Qset}
\blockn{In the previous section we showed that, for any Grothendieck topos $\topos$ endowed with a bound $B$, the quantale $Q=\rel(B)$ already determined the category $\slat(\topos)$. We will now show that if we add\footnote{Actually, because we know that $Q$ is of the form $\rel(B)$, the $^*$ operation is fully determined by the underlying quantale. This comes from the property (Q10) together with this lemma: the condition  $f=g^*$ and $g g^* \leqslant 1$ is equivalent to the condition $\exists u \leqslant 1$, $uf = f$, $gu=g$ , $gf \leqslant 1$ and $u \leqslant fg$ . This lemma is proved using internal logic. } the operation $(\_)^*$ on $Q$, then we can give a complete description of $\topos$ in terms of $Q$.}

\blockn{The theorems \ref{Qset_eq} and \ref{QSet_Grothrepresentation} are the classical results relating Grothendieck topos to Grothendieck quantales they can be found explicitly in  \cite{heymans2010sheaves} and \cite{heymans2012grothendieck} and under different forms in \cite{Freyd_Alegories} and \cite{pitts1988applications}.}

\block{Our starting point will be the following lemmas:

\Lem{Let $X$ be an object of $\topos$, and $Y$ be a sub-quotient of $X$, then the relation on $X$ defined by
 \[ xRy =\text{`` $x$ and $y$ both have an image in $Y$ and these coincide ''} \]
is symmetric ($R^* =R$) and transitive $(R^2 \leqslant R) $. This induces a correspondence between sub-quotients of $X$ and relations $R$ on $X$ such that $R^*=R$ and $R^2 \leqslant R$.
}

\Dem{The symmetry and transitivity of the relation are clear. $Y$ is fully determined by $R$: it is the quotient of $U = \{ x | xRx \}$ by $R$ (which is an equivalence relation on $U$). Conversely, let $R$ be any symmetric and transitive relation on $X$. Let $U = \{ x | xRx \}$, $R$ induces an equivalence relation on $U$, and we have $x R y \Rightarrow x R x$. Hence, $x R y \Leftrightarrow (x \in U) \wedge (x R y) \wedge (y \in U) $, i.e. $R$ is indeed the relation induces by the sub-quotient $U//R$. }}

\block{\Lem{In the situation of the previous lemma, one actually has $R^2 = R$, and the map which sends a sub-object of $Y$ to its pullback in $X$ identifies $\powerob(Y)$ with $R(\powerob(X))$ (where $R$ denotes the endomorphism of sup-lattice of $\powerob(X)$ corresponding to the relation $R$).}

\Dem{Indeed, if $(x R y)$ then $(x R x)$ and $(x R y)$ hence $(x R^2 y)$, this proves that $R \leqslant R^2$, and hence $R = R^2$. Let $P$ be a subset of $X$, $R(P) = \{ x \in X | \exists z \in P, x R z \}$. So $P=R(P)$ if and only if $P$ is included in $U=\{x |x R x \}$ and saturated for the equivalence relation induced by $R$ on $U$. These are exactly the subsets which are pull-backs of subsets of $Y$. }

}

\block{\label{Qset_eq}\Th{The category $\rel(\topos)$ whose objects are the objects of $\topos$ and morphisms from $X$ to $Y$ are sub-objects of $Y \times X$ (the composition being given by the composition of relations) is equivalent to the following category $\proj(Q)$:

\begin{itemize}
\item The objects are the couples $(I,P)$ where $I$ is a set, and $P$ is a matrix in $M_{I,I}(Q)$ such that $P^2 =P$ and $P^* = P$ where $(P^*)_{i,j} = (P_{j,i})^*$.

\item The morphisms from $(J,P')$ to $(I,P)$ are the matrices $M \in M_{I,J}(Q)$ such that $P.M = M$ and $M.P' = M$ (the composition being the product of matrices).

\end{itemize}

Under this equivalence, the opposite of a relation corresponds to the ``trans-conjugation'' of a matrix: $(M^*)_{i,j} = (M_{j,i})^*$.
 }
 
Before proving this theorem We will need one more simple lemma, which is actually the last point of the theorem: 
\Lem{Let $R$ be a sub-object of $(\coprod_{i \in I} B ) \times (\coprod_{j \in J} B)$ corresponding to a morphism $ R: \powerob(B)^I \rightarrow \powerob(B)^J$ represented by a matrix: $(R_{i,j})_{i \in I,j\in J}$, then the opposite relation corresponds to the trans-conjugate matrix $(R^*)_{j,i} = (R_{i,j})^*$ }

\Dem{This can be checked internally: Since $R_{i,j}$ corresponds to the intersection of $R$ with the inclusion $(f_i,f_j)$ of $B \times B$ in $(\coprod_{i \in I} B ) \times (\coprod_{j \in J} B)$, taking the opposite relation will reverse $R_{i,j}$ and exchange the indices. this concludes the proof of the lemma.}
 
\bigskip 
We now prove theorem \ref{Qset_eq}:
\Dem{In order to prove the equivalence of $\proj(Q)$ and $\rel(\topos)$, we will consider a third category $\cat$, the full sub-category of $\slat(\topos)$ of sup-lattices which are of the form $\powerob(X)$ for $X$ an object of $\topos$, and show that both $\proj(Q)$ and $\rel(\topos)$ are equivalent to $\cat$.

\bigskip

The association $X \rightarrow \powerob(X)$ is (one half of) an equivalence from $\rel(\topos)$ to $\cat$. Indeed, it is essentially surjective by definition of $\cat$, and we have already mentioned that morphisms between power objects are the same thing as relations, so it is also fully faithful.

\bigskip

The association $(I,P) \rightarrow P(\powerob(B)^I)$ is (one half) of an equivalence from $\proj(Q)$ to $\cat$. Indeed:

as $B$ is a bound of $\topos$, any object $X$ of $\topos$ is a sub-quotient of some $\coprod_{i \in I} B$, hence there is an endomorphism $F$ of $\powerob(\coprod_{i \in I} B ) = \powerob(B)^I$ such that $F^2 = F$, $F^*=F$, and $\powerob(X) = F(\powerob(B)^I)$. By the result of the previous section, such an endomorphism corresponds exactly to a matrix $P$ such that $(I,P)$ is indeed an object of our category. So this functor is full and well defined (at least on the object). Now a morphism from $P'(\powerob(B)^J)$ to $P(\powerob(B)^I)$ is exactly the data of a matrix $M$ such that $P.M=M$ and $M.P'=M$. This concludes the proof of the equivalences. The last point of the theorem being proved by the lemma.
  
}
}

\block{ \Cor{The topos $\topos$ is equivalent to the (non full) subcategory of $\proj(Q)$, with all objects and with morphisms from $(J,P')$ to $(I,P)$ only the matrices $M$ which satisfy the additional condition: such that $MM^*\leqslant P $ and $P' \leqslant M^*M$.}

\Dem{This two additional conditions indeed characterise functional relations among arbitrary relations, and in a topos functional relations are in correspondence with morphisms.}

}

\block{\label{QSet_Grothrepresentation}\Th{For every Grothendieck quantale $Q$, there exists a topos $\topos$ and a bound $B$ of $\topos$ such that $Q = \rel(B)$.}

Of course, from the previous theorem, such a topos is unique.

\Dem{One could use the construction of $\site(Q)$ given in the introduction, but the proof that this is indeed a site and that it gives back $Q = \rel(B)$ is long and not really illuminating. Instead, we will use results from the theory of allegories (see \cite{Freyd_Alegories}, or \cite{sketches}A.3) which is closely connected to what we are doing here:

In the language of \cite{Freyd_Alegories} a modular quantale $Q$ is a one object locally complete distributive allegory, and $\proj(Q)$ is the systemic completion of $Q$. The result \cite{Freyd_Alegories}$2.434$ proves that $\proj(Q)$ is a power allegory and \cite{Freyd_Alegories}$2.226$ proves that it is has a unit. So in order to apply \cite{Freyd_Alegories}$2.414$ and conclude that $\proj(Q)$ is the category of relations on an elementary topos, we need to proove that it is `tabular'. Using \cite{Freyd_Alegories}$2.16(10)$ it is enough to prove that\footnote{The reader should note that \cite{Freyd_Alegories} uses a reverse composition order for morphisms in category, whereas we use the standard composition order. This explains why the formula we give is different from the one given in the reference.} for each set $X$ the maximal matrix of $M_{X,X}(Q)$ can be written $FG^*$ for  $F,G \in M_{X,Y}(Q)$ with $FF^* \leqslant Id_{Y}$ and $GG^* \leqslant Id_{Y}$. But $(Q10)$ is exactly the assertion that this is true when $X$ is a singleton, and the general case follows easily from $(Q10)$ by taking $Y = X \times I $.

\bigskip

The elementary topos we obtain in this way has arbitrary co-products and is bounded, hence it is a Grothendieck topos.

\bigskip

Finally, if $B$ is the object of $\proj(Q)$ represented by the set $X = \{*\}$ and $P=1$, then $\rel(B)=Q$ and this concludes the proof.
}

}

\block{In the remainder of this section we just give a simpler description of the category $\proj(Q)$ in term of $Q$-Set inspired from the notion of $\mathcal{L}$-sets, when $\mathcal{L}$ is a locale. Our aim is both to provide a formalism suitable for computation and to show that $\proj(Q)$ is exactly a non-commutative generalisation of $\mathcal{L}$-set. We do not know if this formulation has already been presented somewhere or not.

\Def{\begin{itemize}

\item A $Q$-Set is a set $X$ endowed with an application $\Qrel{ \_}{\_}: X \times X \rightarrow Q$ such that:

\[
\begin{array}{c l}
(S1) &   \forall x,y \in X, \Qrel{x}{y} = \Qrel{y}{x}^* .\\
(S2) &\forall x,y,z \in X, \Qrel{x}{y} \Qrel{y}{z} \leqslant \Qrel{x}{z}.
\end{array}
\]

\item A $Q$-relation $R$ from $X$ to $Y$ (two $Q$-sets) is a map: 

\[\begin{array}{c c c}
Y \times X & \rightarrow & Q \\
(y,x) & \mapsto &\Qrel[R]{y}{x} \\
\end{array} \]

such that:

\[ \begin{array}{c l}
(R1) & \Qrel{y}{y'} \Qrel[R]{y'}{x} \leqslant \Qrel[R]{y}{x} \text{ with equality whenever } y = y' \\
(R2) & \Qrel[R]{y}{x'} \Qrel{x'}{x} \leqslant \Qrel[R]{y}{x} \text{ with equality whenever } x= x'.
\end{array}
\]

\item A $Q$-function from $X$ to $Y$ is a $Q$-relation:

\[\begin{array}{c c c}
Y \times X & \rightarrow & Q \\
(y,x) & \mapsto &\Qfun{y}{f}{x} \\
\end{array} \]

which (in addition to $(R1)$ and $(R2)$) satisfies:

\[ \begin{array}{c l}
(F1) &\Qfun{y}{f}{x} \Qfun{y'}{f}{x}^* \leqslant \Qrel{y}{y'} \\
(F2) & \Qrel{x}{x} \leqslant \bigvee_{y} \Qfun{y}{f}{x}^* \Qfun{y}{f}{x}
\end{array}
\]

\item $Q$-relations and $Q$-functions can be composed by the formula:
\[ \Qrel[RQ]{z}{x} = \bigvee_y \Qrel[R]{z}{y} \Qrel[Q]{y}{x} \]
\[ \Qfun{z}{f \circ g}{x} = \bigvee_{y} \Qfun{z}{f}{y} \Qfun{y}{g}{x}\] 

\item The opposite of a $Q$-relation is given by:

\[ \Qrel[R^*]{x}{y} = \Qrel[R]{y}{x}^* \]
\end{itemize}
}
}

\block{
\Prop{Consider the following modification of the axioms:

\begin{itemize}
\item[(S2')] $\Qrel{x}{y} = \bigvee_{t} \Qrel{x}{t} \Qrel{t}{y}$.
\item[(R1')] $\bigvee_{y'} \Qrel{y}{y'} \Qrel[R]{y'}{x} = \Qrel[R]{y}{x}$ 
\item[(R2')] $\bigvee_{x'} \Qrel[R]{y}{x'} \Qrel{x'}{x} = \Qrel[R]{y}{x}$
\item[(F2')] $\Qrel{x}{x'} \leqslant \bigvee_{y} \Qfun{y}{f}{x}^* \Qfun{y}{f}{x'}$
\end{itemize}

Then assuming $(S1)$ holds, $(S2)$ and $(S2')$ are equivalent.
 
And assuming $X$ and $Y$ are $Q$-sets, $(R1)$ is equivalent to $(R1')$, $(R2)$ is equivalent to $(R2')$ and assuming additionally $(R2)$ then $(F2)$ is equivalent to $(F2')$.

\bigskip

In particular $Q$-Sets are exactly the same as objects of $\proj(Q)$, and $Q$-relations and $Q$-functions correspond respectively to morphisms in $\proj(Q)$, and morphisms which are sent to functional relations by the equivalence of \ref{Qset_eq}.
}

We will need the following lemma:

\Lem{In any $Q$-set, one has
\[ \Qrel{x}{y} \Qrel{y}{y} = \Qrel{x}{y} \]
\[ \Qrel{x}{x} \Qrel{x}{y} = \Qrel{x}{y}\]
}
\Dem{Indeed (for the second equality), for any $q \in Q$ element of a modular quantale one has:

\[ q \leqslant (1.q \wedge q) \leqslant (1 \wedge qq^*)q \leqslant qq^*q \]

So:
\[ \Qrel{x}{y} \leqslant \Qrel{x}{y} \Qrel{y}{x} \Qrel {x}{y} \leqslant \Qrel{x}{x} \Qrel {x}{y} \]

The reverse inequality being a consequence of $(S2)$, one has the desired equality.}

We now prove the theorem:

\Dem{
\begin{itemize}
\item Clearly, $(S2')$ implies $(S2)$. Assume that $(S2)$ and $(S1)$ hold, hence that $X$ is a $Q$-set. One can apply the lemma and one has:

\[ \Qrel{x}{y} \leqslant \bigvee \Qrel{x}{t}\Qrel{t}{y} \]

by taking $t=x$ or $t=y$. The reverse inequality follows from $(S2)$.

\item $(R1) \Rightarrow (R1')$ is clear because of the equality case. Assuming $(R1')$ one has immediately $ \Qrel{y}{y'} \Qrel[R]{y'}{x} \leqslant \Qrel[R]{y}{x} $. So we just have to prove that  $ \Qrel{y}{y} \Qrel[R]{y}{x} = \Qrel[R]{y}{x} $. But:

\[ \Qrel{y}{y} \Qrel[R]{y}{x} = \bigvee_{y'} \Qrel{y}{y} \Qrel{y}{y'} \Qrel[R]{y'}{x} = \bigvee_{y'} \Qrel{y}{y'} \Qrel[R]{y'}{x} = \Qrel[R]{y}{x} \]

The equivalence of $(R2)$ and $(R2')$ is proved the same way.

\item $(F2)$ is a special case of $(F2')$. Assume $(F2)$ then:

\[ \begin{array}{r l}
\Qrel{x}{x'} = \Qrel{x}{x} \Qrel{x}{x'} & \leqslant \bigvee_y \Qfun{y}{f}{x}^* \Qfun{y}{f}{x} \Qrel{x}{x'} \\ &\leqslant \bigvee_y \Qfun{y}{f}{x}^* \Qfun{y}{f}{x'}
\end{array} \]

\end{itemize}

The fact that $Q$-Sets are exactly the same as objects of $\proj(Q)$, and $Q$-relations and $Q$-functions correspond respectively to morphisms in $\proj(Q)$, and morphisms which are sent to functional relations by the equivalence of \ref{Qset_eq} is now immediate: If we replace the original axioms by these modified version, and if we interpret $\Qrel{x}{y}$,$\Qrel[R]{x}{y}$ and $\Qfun{x}{f}{y}$ as matrix coefficients then the conditions imposed on them are exactly those for being objects and morphisms of $\proj(Q)$.

}

}

\subsection{Relational representations of Grothendieck quantales}
\label{RelRep}
\blockn{In the previous section we constructed a topos $Q$-Sets from a Grothendieck quantale $Q$. In this section we describe the theory classified by this topos, that is study the morphisms from an arbitrary topos $\topos$ to the topos of $Q$-sets. This also explains in which sense the equivalence between Grothendieck quantales and Grothendieck toposes is functorial.}

\block{\label{RelRep_morphism}\Def{A morphism of modular quantales is an application $f:Q \rightarrow Q'$ between two modular quantales such that:

\begin{itemize}
\item $f$ commutes to arbitrary supremum (in particular it preserves the smallest elements)
\item $f$ commutes to finite intersections (in particular it preserves the top element $\top$).
\item $f$ is a morphism of unitary monoids (in particular it preserves $1$).
\item $f$ commutes to the involution.
\end{itemize}

A Relational representation of a modular quantale $Q$ is the datum of an inhabited set $X$ endowed with a modular quantale morphism $\pi$ from $Q$ to $\rel(X)$. A morphism of relational representations is a map from $X$ to $X'$ such that for each $q \in Q$ if $(x,y) \in \pi(q)$ then $(f(x),f(y)) \in \pi'(q)$ .

}

}

\renewcommand{\toposdefaut}{1}
\block{\Th{The topos of $Q$-sets classifies the relational representations of $Q$, the universal representation being given by the action of $Q$ on the bound $B$ (which corresponds to the $Q$-set $\{*\}$ with $\Qrel{*}{*} = 1$).
In other words, if $\topos$ is any topos then there is an equivalence of categories between the geometric morphisms from $\topos$ to $Q$-sets, and the relational representations of $Q$ inside $\topos$.
And this equivalence is given by $f \mapsto f^*(B)$.} 

This theorem is essentially the same as   theorem 2.9 of \cite{pitts1988applications}.

\Dem{As a bound, the object $B$ has to be in particular inhabited, hence it is indeed a relational representation. So any geometric morphism from $\topos$ to $Q$-sets does yield a relational representation of $Q$ on $f^*(B)$ and any natural transformation gives a morphism of representation. So the functor mentioned in the theorem indeed exists.

\bigskip

If $f$ is a geometric morphism from $\topos[1]$ to $Q$-sets, then $f^*$ induces a $\slat$-enriched functor from $\proj(Q)$ to $\rel(\topos)$.

Because $\proj(Q)$ is generated by $B$ under co-product and splitting of projection, and since by \ref{Slcat_biproduct} arbitrary co-products (as well as spliting of projection) are universal co-limits in $\slat$-enriched category, any relational representation $(X,\pi)$ of $Q$ in a topos $\topos$ extends in a uniquely defined $\slat$-enriched functor from $\proj(Q)$ to $\rel(\topos)$: one havehas to send the couple $(I,P)$ on $\pi(P) \coprod_{i \in I} X$, and any morphism in $\proj(Q)$ is a matrix which is sent to the matrices ``$\pi(M)$'' defining a relation in $\topos$.

Moreover if $f$ and $g$ are two geometric morphisms from $\topos$ to $Q$-sets, then morphisms between the relational representation they induce uniquely extend to natural transformations between $f^*$ and $g^*$.

So we just have to prove that if $(X,\pi)$ is a relational representation of $Q$, then the induced functor $v$ from $\proj(Q)$ to $\rel(\topos)$ comes from a geometric morphism from $\topos$ to $Q$-sets.

\begin{itemize}

\item As $\pi$ commutes with $*$, so does $v$. Hence $v$ preserves functional relations and induces a functor from $Q$-sets to $\topos$.

\item The terminal object of $Q$-sets is the quotient of $B$ by its maximal relation, and since $\pi$ preserves the maximal relation, the terminal object of $Q$-sets is sent to the quotient of $X$ by its maximal relation, which is the terminal object of $\topos$ because $X$ is inhabited. So $v$ preserves the terminal object.

\item Let 

\[
\begin{tikzcd}[ampersand replacement=\&]
P \arrow{r} \arrow{d} \& X \arrow{d}{f} \\
Y \arrow{r}{g} \& S \\
\end{tikzcd}
\]

 be a pull back diagram in $Q$-sets, then $P$ can be identified with the relation $f^* g$ on $X \times Y$, indeed internally $f^* g$ is the relation $\{(x,y | f(x)=g(y) \}$ hence it is the fiber product $X \times_S Y$.
So $v$ preserves pull-back. As $v$ preserves the terminal object, it preserves all limit.

\item In a topos, a collection of maps $f_i: A_i \rightarrow A$ is a covering if and only if $1_A \leqslant \bigvee_{i \in I} f_i f_i^*$ as $v$ preserves all the structures involved, $v$ preserves covering families.

\end{itemize}

All these properties together imply that $v$ is indeed the $f^*$ functor of a geometric morphism and conclude the proof.

}

}
\renewcommand{\toposdefaut}{0}
\subsection{Internally bi-linear maps between $Q$-modules}
\label{SlatBil}
\blockn{The category $\slat(\topos)$ is endowed with a tensor product. In the case where $\topos = \sh(\mathcal{L})$ is the topos of sheaves on a frame $\mathcal{L}$, one can see that through the identification of $\slat(\topos)$ with $\module_\mathcal{L}$ this tensor product corresponds to the natural tensor product over $\mathcal{L}$, which as in the case of commutative algebras is defined by the universal property: the maps from $M \otimes_\mathcal{L} N$ to $P$ are the bi-linear morphisms from $M \times L$ to $P$ such that for all $l \in \mathcal{L}$, $f(m,n.l)=f(m.l,n)=f(m,n).l$.

The main result of this section is that, in the general case, even if the tensor product of two $Q$-modules can be difficult to compute explicitly, the set $\bil_{\topos}(M\times N, P)$ of internal bilinear map from $M \times N$ to $P$ has a strikingly simple description in term of the corresponding right $Q$-modules. This leads in particular to a simple description of the category of internal locales of $\topos$ in terms of a Grothendieck quantale representing $\topos$. More precisely:}

\block{\label{SlatBil_Def}\Def{If $A$, $B$ and $C$ are three right modules over a Grothendieck quantale $Q$, we say that a map $f:A \times B \rightarrow C$ is $Q$-bilinear if it is a bi-linear morphism of sup-lattices and if it satisfies the following three conditions:

\begin{enumerate}
\item $\tilde{f}(aq,b) \leqslant \tilde{f}(a,bq^*)q$
\item $\tilde{f}(a,bq) \leqslant \tilde{f}(aq^*,b)q$
\item $\tilde{f}(a,b).q \leqslant \tilde{f}(aq,bq)$.
\end{enumerate}

We will denote by $\bil_Q(A \times B,C)$ the set of $Q$-bilinear maps.
}

$\bil_Q(A \times B,C)$ is a sup-lattice for the pointwise ordering (with supremum computed pointwise), and it is an $\slat$-enriched functor in each of the three variables (contravariant in the first variables) with the functoriality given by composition.

The main result of this section (theorem \ref{SlatBil_theorem}) is that this functor is isomorphic to the functor of internal bilinear maps.

}

\block{
\label{SlatBil_construction}Let $M$,$N$ and $P$ be internal sup-lattices in $\topos$, let $\widetilde{M}$, $\widetilde{N}$ and  $\widetilde{P}$ be the corresponding right $Q$-modules (i.e. $\tilde{M} = \hom_{\topos}(B,M)$).
Let $f$ be a bilinear morphism from $M \times N$ to $P$.

Then one can define a map $\tilde{f}$ from $\widetilde{M} \times \widetilde{N}$ to $\widetilde{P}$ by the (internal) formula:

\[ \tilde{f}(m,n):= b \mapsto f(m(b),n(b)) \in P \]

With $m$ and $n$ elements of $\widetilde{M}$ and $\widetilde{N}$, that is, maps from $B$ to $M$ and $N$, then $\tilde{f}(m,n)$ is indeed an element of $\widetilde{P} = \hom_{\topos}(B,P)$. 

\Prop{The map $\tilde{f}$ is a $Q$-bilinear morphism map in the sense of \ref{SlatBil_Def}. Moreover the construction $f \rightarrow \tilde{f}$ defines a morphism of $\slat$-enriched functors:

\[ \mu(M,N,P) : \bil_{\topos}( M \times N,P) \rightarrow \bil_Q(\widetilde{M} \times \widetilde{N}, \widetilde{P}) \]

 }

\Dem{The (sup-lattice) bilinearity is immediate: supremum in $\widetilde{M},\widetilde{N}$ and $\widetilde{P}$ corresponds to pointwise internal supremum hence the bilinearity of $\tilde{f}$ simply comes from the internal bilinearity of $f$.

Recall that by definition one has internally for any $m \in \widetilde{M}$ and $q\in Q=\rel(B)$: $\displaystyle m.q(b) = \bigvee_{(b',b)\in q} m(b')$. 

All three properties defining $Q$-bilinearity are then easily checked internally:

\begin{enumerate}
\item \[ \tilde{f}(mq,n)(b) = \bigvee_{(b',b)\in q} f( m(b'),n(b)) \]
Whereas:

\[ \tilde{f}(m,nq^*)q(b) = \bigvee_{(b',b)\in q} \tilde{f}(m,nq^*)(b') = \bigvee_{(b',b) \in q, (b',b'')\ in q } f(m(b'),n(b'')) \]

So the first term corresponds to the restriction of the union to $b=b''$ of the second and is indeed smaller.

\item Same proof.

\item \[ \left[ \tilde{f}(m,n)q \right](b) = \bigvee_{(b',b)\in q }f(m(b'),n(b')) \]

Whereas 

\[ \tilde{f}(mq,nq)(b) = f(mq(b),nq(b))=\bigvee_{(b',b) \in q, (b'',b)\in q } f(m(b''),n(b')) \]

So the first term corresponds to the restriction of the union to $b'=b''$ of the second and is indeed smaller.

\end{enumerate}

Also $f \mapsto \tilde{f}$ commutes to supremum, because if one takes $f_i$ a (external) family of internal bilinear maps, $m \in \widetilde{M}$ and $n \in \widetilde{N}$ then (internally) for any $b \in B$:

\[ \left(\widetilde{\bigvee_i f_i } \right)(m,n)(b) = \bigvee_i f_i(m(b),n(b)) = \left(\bigvee_i \tilde{f_i}(m,n) \right) (b) \]

And the functoriality is immediate: $\tilde{f}(m,g(n)) := b \mapsto f(m(b),g(n(b))$ is indeed the map attached to $f( \_,g(\_))$ and $\tilde{g}(\tilde{f}(m,n) := b \mapsto g(f(m(b),n(b)) = \widetilde{g \circ f}(m,n)$.

}

}

\block{\label{SlatBil_theorem}\Th{The construction $f \mapsto \tilde{f}$ from \ref{SlatBil_construction}, defines an isomorphism of $\slat$ enriched functors:

\[  \mu : \bil_{\topos}(M \times N,P) \simeq \bil_Q(\widetilde{M} \times \widetilde{N},\widetilde{P}). \]  }

The functoriality of the association and the fact that it commutes to supremum have already been mentioned, so it only remains to prove that it is a bijection.
The proof of this theorem will be completed in \ref{SlatBil_proof} after proving a few lemmas.

}

\block{\label{SlatBil_injective}\Lem{The association $f \mapsto \tilde{f}$ of \ref{SlatBil_construction} is injective.} 

\Dem{Let $f$ and $g$ be two internal bi-linear maps from $M \times N$ to $P$ such that $\tilde{f} = \tilde{g}$. This means that for each map $(m,n): B \rightarrow M \times N$ one has internally:

\[\forall b \in B, f(m(b),n(b))= g(m(b),n(b))\]

i.e.:

\[ f \circ (m,n) = g \circ (m,n). \]

But we already explained in the last part of the proof of \ref{Slcat_topos} that any map from a sub-object $U$ of $B$ to a sup-lattice can be extended (canonically) to a map on all of $B$. As $B$ is a bound, maps from sub-objects of $B$ can cover $M \times N$ and by the extension arguments, maps from $B$ cover $M \times N$, so we can conclude from the previous formula that $f=g$.

}
}

\block{\label{SlatBil_lemma2}\Lem{Let $h: Q \times Q \rightarrow \tilde{P} \in \bil_Q(Q \times Q, \tilde{P} )$ where $Q$ is endowed with its right action on itself.
Then:
\begin{itemize}
\item Let $c\in \tilde{P}$ and $a,b \in Q$ such that $a.a^* \leqslant 1$, $b.b^* \leqslant 1$. If one has $c \leqslant h(a,b)$ then for all $x,y \in Q$:

\[ c.(a^*.x \wedge b^*.y) \leqslant h(x,y) \]

\item For all $x,y$ one has:

\[ h(x,y) = \bigvee_{aa^* \leqslant 1, bb^* \leqslant 1} h(a,b)(a^*x \wedge b^* y) \]

\end{itemize}
 }

\Dem{
For the first point:

Let $t = (a^*x \wedge b^*.y) \in Q$. Then one has:

\[a.t \leqslant (aa^* x \wedge ab^* y) \leqslant aa^* x \leqslant x \]
\[b.t \leqslant (ba^* x \wedge bb^* y) \leqslant bb^* y \leqslant y \]

So:
\[ c.t \leqslant h(a,b).t \leqslant h(a.t,b.t) \leqslant h(x,y) \]

For the second point:

Let $h'(x,y) = \bigvee_{aa^* \leqslant 1, bb^* \leqslant 1} h(a,b)(a^*x \wedge b^* y)$.

Clearly, $h'$ is also in $\bil_Q(Q \times Q, \tilde{P})$.

The first point shows that $h' \leqslant h$. For the reverse inequality we will proceed in several steps:

\begin{itemize}

\item If $(xx^* \leqslant 1) , (yy^* \leqslant 1)$. Let $D(x) = 1 \wedge x^* x$ and $D(y)= 1\wedge y^* y$. We note that for elements smaller than $1$, the involution is the identity and composition and intersection coincide (these can be proved by applying the modularity law, or by using theorem \ref{QSet_Grothrepresentation} as a black box and checking it internally in the corresponding topos). So

\[x.D(x) = x(1 \wedge x^* x) \geqslant (x \wedge x) = x \]

hence $x.D(x) = x$ and $y.D(y)=y$ also,

\[ x^* x \wedge y^* y \geqslant D(x) \wedge D(y) = D(x).D(y) \]

Aslo for any $e \leqslant 1$:

\[h(xe,y) \leqslant h(x,ye)e \leqslant h(x,y)e \leqslant h(xe,ye) \leqslant h(xe,y)  \]

hence
\[ h(xe,y)=h(x,ye)=h(x,y)e \]

Finally:

\[ h(x,y)(x^*x \wedge y^*y) \geqslant h(x,y)D(x)D(y) = h(xD(x),yD(y))= h(x,y) \]

So $h'(x,y) \geqslant h(x,y)$

\item We will now assume that $x$ is arbitrary and $y$ is simple. As $Q$ is a Grothendieck quantale, $x$ can be written as a supremum of elements of the form $uv^*$ with $u$ and $v$ simple. So, by bi-linearity of $h$, it is enough to prove that $h(x,y) \leqslant h'(x,y)$ when $x$ is of the form $uv^*$. In this case:

\[ h(uv^*,y) \leqslant h(u,yv)v^* = h'(u,yv)v^* \leqslant h'(uv^*,yvv^*) \leqslant h'(uv^*,y) \] 

\item If both $x$ and $y$ are now arbitrary, then the same technique allows ones to conclude.

\end{itemize}

 }

}

\block{\label{SlatBil_firstcase}\Cor{Theorem \ref{SlatBil_theorem} holds whenever $M = N =\powerob{B}$ (ie. $\mu(\powerob(B),\powerob(B),\tilde{T})$ is an isomorphism)}
\Dem{Injectivity is known by lemma \ref{SlatBil_injective}. So we just have to prove surjectivity. Let $h$ be any bi-linear map from $Q \times Q = \tilde{M} \times \tilde{N} \rightarrow \tilde{P}$ satisfying our three conditions. Then, by lemma \ref{SlatBil_lemma2}, one has that:

\[ h(x,y) = \bigvee_{aa^* \leqslant 1, bb^* \leqslant 1} h(a,b)(a^*x \wedge b^* y). \]

One can see that maps of the form $(x,y) \mapsto t(ux \wedge vy)$ with $u,v \in Q$ and $t\in \tilde{T}$ correspond to the internal bi-linear map which sends $(p,q) \in \powerob(B) \times \powerob(B)$ to $t(u(p) \wedge v(q)$ where $u$ and $v$ are seen as endomorphisms of $\powerob(B)$ and $t$ as a morphism from $\powerob(B)$ to $T$. Indeed if $g(p,q)=t(u(p) \wedge v(q))$ then

\[ \tilde{g}(x,y) = t(ux(b) \wedge vy(b)) = \bigvee_{(b',b) \in (ux \wedge vy)} t(b') = t(ux \wedge vy) \]

Hence $h$ can be written as a supremum of maps coming from internal bi-linear maps, but as $f\mapsto \tilde{f}$ commutes to arbitrary supremum, this shows that $h$ does also come from a bilinear map.
}
}

\block{\label{SlatBill_extension}At this point, we have two possibilities. We can conclude by an argument of extension by inductive limit, or use the following argument which we found more convincing: 

\Prop{Assume that $\mu(\powerob(B),M,P)$ is an isomorphism for some $M$ and $P$ in $\slat(\topos)$. Then $\mu(N,M,P)$ is an isomorphism for any $N$ in $\slat(\topos)$. }

\Dem{We already know that $\mu$ is injective. Hence it remains to show the surjectivity.

Let $M$ and $P$ such that $\mu(\powerob(B),M,P)$ is an isomorphism. the internal sup-lattice $[M,P]$ then corresponds to the right $Q$-module:
\[\begin{array}{r c l} \hom_{\topos}(B,[M,P])& =& \hom_{\slat(\topos)}(\powerob(B),[M,P]) = \bil_{\topos}(\powerob(B) \times M,P) \\ &= &\bil_Q(Q \times \tilde{M}, \tilde{P}) \end{array}\]

where the action of $Q$ on the last term is given by the left action of $Q$ on itself.
 
Now let $g \in \bil_Q(\tilde{N} \times \tilde{M},\tilde{P})$. For any $n\in \tilde{N}$, the map:

 \[ g_n : (q,m) \mapsto g(nq,m) \]
 
is an element of $\bil_Q(Q \times \tilde{M}, \tilde{P})$. The map $(n \mapsto g_n)$ is a morphism of right $Q$-modules and hence internally corresponds to a map $f:N \rightarrow [M,P]$ which in turn corresponds to a map $f \in \bil_{\topos}(N \times M,P)$. Finally $\tilde{f}=g$ because for any $n \in \tilde{N}$, $(q,m) \mapsto \tilde{f}(nq,m)$ is by construction of $\tilde{f}$ the map $g_n \in \bil_Q(Q \times \tilde{M},\tilde{P})$ and hence $\tilde{f}$ agrees with $g$. This concludes the proof.

}

}

\block{\label{SlatBil_proof}We can now finish the proof of \ref{SlatBil_theorem}: by \ref{SlatBil_firstcase} we know that for any $T \in \slat(\topos)$, $\mu(\powerob(B), \powerob(B),T)$ is an isomorphism. Hence by \ref{SlatBill_extension}, $\mu(N, \powerob(B),T)$ is an isomorphism for each $N$ and $T$, as one can freely exchange the first two variables, $\mu( \powerob(B),N,T)$ is also an isomorphism, and a second application of \ref{SlatBill_extension} allows one to conclude.
}

\block{\Cor{\label{SlatBil_locale}The category of internal locales of $\topos$ (equivalently, the category of toposes which are localic over $\topos$) is equivalent to the category of locales $\mathcal{L}$ endowed with a right action of $Q$ such that:
\begin{itemize}
\item As a sup-lattice $\mathcal{L}$ is a right $Q$-module.
\item One has the modularity condition:
\[\forall m,n \in \mathcal{L}, \forall q \in Q, m\wedge nq \leqslant (mq^* \wedge n)q \]

\end{itemize}

We will call such an action a modular action.
(the morphisms of this category being the $Q$-equivariant morphisms of locale).
} 

Of course, the equivalence is given by the usual functor $\mathcal{L} \mapsto \hom_{\topos}(B,\mathcal{L})$.

\Dem{Let $\mathcal{L}$ be a locale in $\topos$ then $\tilde{\mathcal{L}}$ is indeed endowed with an operation $\tilde{\wedge}$ which is $Q$-bilinear. Also:

\[ \tilde{\wedge}(m,n)(b) = m(b) \wedge n(b) \]

is the intersection of $\hom_{\topos}( B, \mathcal{L})$, hence $\mathcal{L}$ is indeed a locale and intersection is indeed $Q$-bilinear.

\bigskip

Conversely, if $\mathcal{L}$ is a locale endowed with a modular action of $Q$, then the operation $\wedge$ is $Q$-bilinear (the second axiom comes from the symmetry, and the third axiom because multiplication by $q$ is order preserving). Hence $\mathcal{L}$ corresponds to an internal sup-lattice $L$ equipped with a bi-linear map $m$ coming from $\wedge$. This bi-linear map has to be the intersection map of $L$ because both $m$ and $\wedge_L$ induce the same map when we externalise it by watching the morphisms from $B$ to $L \times L$, and the proof of the injectivity of the externalisation process done in \ref{SlatBil_injective} works without assuming that the map $f$ is bilinear.  

}

}

\block{Also, from the description of $\bil_{\topos}(M \times N ,P)$ it is possible to obtain an explicit description of both the tensor product and the internal hom objects in terms of the corresponding $Q$-modules: for the hom-object it has been done in the proof of \ref{SlatBill_extension} and for the tensor product, the description of $\bil_{\topos}(M \times N ,P)$ translates into a presentation by generators and relations of $\widetilde{M \otimes N}$. In order to completely handle the monoidal structure of $\slat(\topos)$ in terms of the category of $Q$-modules it remains to understand what $\widetilde{M^*}$ is.

An element of $\widetilde{M^*}$ is an application from $B$ to $M^*$, hence it is the same thing as an application from $B$ to $M$ (i.e. an element of $\widetilde{M}$) but with the reverse order relation, hence as a sup-lattice $\widetilde{M^*}=\widetilde{M}^*$. A simple computation shows that $q \in Q$ acts on $\widetilde{M}^*$ by the adjoint of the action of $q^*$ on $M$.

}

\blockn{We also note that the notion of $Q$-bilinear map makes sense when $Q$ is only a modular quantale, that one can define a ``tensor product'' $M\otimes_Q N$ universal for $Q$-bilinear maps from $M\times N$. But in general this tensor product fails to be associative. }

\subsection{Representations of modular quantales}
\block{\label{ModQ_conj}Conjecture: Every modular quantale is of the form $\rel(X)$ with $X$ an object of a topos. }

\block{Let $Q$ be any modular quantale, we can try to consider the classifying topos of the theory of relational representations of $Q$ and hope that the quantale of relations on the universal representation is isomorphic to $Q$. Unfortunately this is not true in general as the following example shows:

\bigskip

Let $n$ be any integer $>2$, and $G=S_n$ be the permutation group and consider the set $X_n = \{1,\dots,n\}$ endowed with its natural action of $S_n$ as a $S_n$-set.
The only relations on $X_n$ (in $S_n$-sets) are: $\emptyset, 1, \top, \Delta$ where $\Delta$ is the complementary relation of $1$.

As $n >2$ one has $\Delta^2=\top$. Let $Q$ be the quantale of relations on $X$ in $S_n$-sets. $Q$ does not depend on $n$. A relational representation of $Q$ (in a topos) is just an object $S$ which has a decidable equality (the diagonal sub-object is complemented) and at least three distinct elements. The universal model $U$ of this theory has more relations on it than just the four elements of $Q$:

Indeed, let

 \[ P = \{(x,y) \in U| \exists x_1,\dots, x_4 \in U^4 \text{ pairwise disjoint in U} \}.\]

then the pull-back of $P$ in the representation $X_3$ is $\emptyset$ whereas the pull-back of $P$ in the representation $X_n$ for $n>3$ is $\top$ hence $P$ cannot be any of the objects of $Q$ in the universal representation.

}

\subsection{Towards a convolution $C^*$ algebra attached to a quantale}

\label{ConvProd}

\blockn{ In the special case where the topos $\topos$ is an \'etendue, ie the topos of equivariant sheaves over an \'etale (localic) groupoid $G =(G_0,G_1)$, the quantale associated to the bound $B$ such that $T_B = G_0$ is the set of open subsets of $G_1$, and the composition law is given by the direct image of open subsets in the composition law of $G_1$.

\bigskip 

With this fact and the usual construction of a $C^*$-algebra from a groupoid (see \cite{paterson1999groupoids} ) in mind it is natural to try to construct a $C^*$-algebra from a quantale by defining a convolution product over a subset of continuous functions ``over $Q$'' (or over open subspaces of $Q$). Indeed we can view $Q$ as a locale by forgetting its composition law and involution, then use the involution to get an involution on continuous functions and hope\footnote{This will not be the case in full generality.} that the composition law on $Q$ will allow us to construct a convolution product on continuous functions. 
}

\blockn{In order to perform this construction the general idea is the following: a continuous function on $Q$ is the same thing as a function from $B \times B$ to the sheaf of real or complex numbers. Hence it can be thought of, internally, as an infinite matrix whose rows and columns are indexed by $B$, and we can use the multiplication of matrices to define the convolution product. If the coefficients are all positive and we allow infinite coefficients the product should be always defined. There are two difficulties that arise when we try to define the matrix product internally:

\begin{itemize}

\item Matrix multiplication requires a summation indexed by the elements of $B$. It can not be done if we don't assume that $B$ has a decidable equality. The reason for this is that without this assumption, when we look at partial sums $f(b_1) + \dots +f(b_n)$, for $(b_1, \cdots b_n)$ elements of $B$ we cannot say whether the elements $b_i$ are distinct, hence we cannot assert that we have not counted some value of $f$ twice. 

\item The sum of an infinite number of terms will in full generality be defined as the supremum of all the possible finite sums. In general the object $\Rt$ of continuous real numbers (i.e. two sided Dedekind cut) does not always have supremum. In order to define a supremum we need to replace the usual ``continuous'' real numbers by the lower semi-continuous real numbers (one-sided Dedekind cut), so the result of the convolution will in general be a lower semi-continuous function.

\end{itemize}
}

\block{We now move to the precise definition: 

\Prop{Let $B$ be an object of a topos $\topos$ with a decidable equality. Let (internally) $f$ and $g$ be functions from $B \times B$ to $\Rlscp$. Then we can define a function $(f * g)$ from $B \times B$ to $\Rlscp$ such that for all $q \in \Qt$,  one has 
\[ q < (f*g)(b,b') \]

if and only if: $ \exists n \in \NN, b_1, \cdots ,b_n \in B$ such that:
\[ \forall i \neq j, B_i \neq B_j\]
and
\[q < \sum_{i=1}^n f(b,b_i) g(b_i,b') \]

where $\displaystyle q < \sum_{i=1}^n f(b,b_i) g(b_i,b')$ naturally means:

$\exists u_1,\cdots,u_n, v_1,\cdots v_n \in \Qt$,  such that $ u_i < f(b,n_1), v_i < g(b_i,b')$  and  $q <  \sum_{i=1}^n u_i v_i$.
}

\Dem{ All we have to do is check that the set $X$ of $q$ such that:

\[ \exists n \in \NN, b_1, \cdots ,b_n \in B, \text{ such that } \forall i \neq j, B_i \neq B_j \text{ and } q < \sum_{i=1}^n f(b,b_i) g(b_i,b') \]

is indeed a positive one-sided Dedekind cut.

It is positive because, by taking $n=0$ all negative $q$ are in $X$. If $q'<q$ and $q \in X$ then clearly $q' \in X$. If $q \in X$, then $q < \sum f(b,b_i) g(b_i,b')$ so there exists $q' >q$ such that $q' < \sum f(b,b_i) g(b_i,b')$ hence $q' \in X$. This concludes the proof.

}

}

\block{\Prop{The convolution product defined is associative and the characteristic function of the unit of $Q$ is a unit.}

\Dem{This an immediate consequence of the fact that internally, the composition of matrices is associative and that the identity matrix is a unit for the composition of matrices. All we need  is a constructive version of Fubini's theorem for sums (indexed by decidable sets) of positive lower semi-continuous real numbers. The usual proof can easily been made constructive, or we can apply the general Fubini's theorem proved in \cite{vickers2011monad}.  } 
 }
 
\blockn{This is interpreted externally as the construction of a convolution product on the set of lower semi-continuous functions on the underlying space of $Q$. This corresponds exactly to the construction of the convolution algebra of an étale groupoid (see for example \cite{paterson1999groupoids}\footnote{In this reference, étale groupoids are called r-discrete groupoids.} for this construction).

It should be possible to obtain something similar to the more general construction of the convolution algebra of a locally compact groupoid endowed with a Haar system, by replacing the bound $B$ by an internal space endowed with an internal ``measure'' (a valuation to be more precise). This would require using internal measure theory as developed in \cite{vickers2008localic} and \cite{vickers2011monad}.} 

\blockn{Of course, in full generality this kind of construction may fail to give anything interesting: it may not restrict to an operation on continuous functions, or it may even yield an everywhere infinite result. The last section will give a more precise picture of the situation in the special case where the underlying topological space of $Q$ is discrete.  }

\section{Atomic Quantales}

\subsection{Introduction}

\blockn{In this section we focus on a really special case of the theory explained in the previous section: when the underlying space of the quantale $Q$ (that is, the space obtained by forgetting the product on $Q$ and seeing it as a locale) is a discrete topological space. This corresponds to the study of atomic toposes (see \cite{sketches} C.3.5 for the theory of atomic toposes).

In this situation, the convolution product constructed at the end of the previous section is easier to understand, and among other things we will explain a simple necessary and sufficient condition for this convolution product to be defined and interesting on continuous functions. In this special case two additional features appear:  a canonical `time evolution' on the $C^*$ algebra obtained this way, and we will observe that when the convolution product is well defined it restricts into a product on finitely supported integer valued functions, which gives rise to an algebra over $\mathbb{Z}$ with possibly interesting arithmetic and combinatorial properties. This algebra can be interpreted as a subalgebra of the algebra of endomorphisms of the free $\Zt$ modules of base $B$.}

\blockn{All the constructions will be made constructively, but at some point we will need to take the assumption that the underlying space of $Q$ is decidable. }

\subsection{Atomic quantales and atomic toposes}
\label{At}
\blockn{An object $X$ of a topos is said to be an atom if $\sub(X) \simeq \soc[]$\footnote{Assuming classical logic in the base topos, this means that $X$ is non-empty and has no non trivial sub-object.}. A direct image of an atom by a morphism is again an atom. A topos is said to be atomic if it satisfies one of the following equivalent properties:

\begin{itemize}
\item The atoms form a generating family.
\item For every object $X\in |\topos|$, $\sub(X)$ is an atomic locale (i.e. of the form $\powerob(S)$ form some $S$.
\item $\topos$ is the topos of sheaves over an atomic site (i.e. a site such that the covering sieves are exactly the inhabited sieves).
\end{itemize}

For details on the theory of atomic toposes one can consult \cite{sketches} C3.5\footnote{This reference studies atomic geometric morphisms, but a geometric morphism $f:\topos[1] \rightarrow \topos$ is atomic if and only if $\topos[1]$ is atomic as a $\topos$-topos. }.
}

\block{\Prop{Let $\topos$ be a topos and $B$ be a bound of $\topos$, then the following conditions are equivalent:

\begin{enumerate}

\item $\topos$ is atomic over sets. 

\item $Q = \rel(B)$ is (as a locale) atomic, i.e. its underlying poset is of the form $\powerob(X)$ for some set $X$.

\item $Z(Q) = \sub(B)= \{ q \in Q |  \leqslant 1 \}$ is atomic.

\end{enumerate}

}

\Dem{ $1. \Rightarrow 2.$ is clear, because in an atomic topos all the lattices $\sub(X)$ are atomic. The implication $2. \Rightarrow 3.$ is also clear.
We will prove $3. \Rightarrow 1.$: Every object of $\topos$ can be covered by subobjects of $B$, so if every subobject of $B$ can be covered by atoms (which is the assumption in $3.$) every object of $\topos$ can be covered by atoms and hence $\topos$ is atomic. This concludes the proof.}

Of course, if $\topos$ is an atomic topos and $X$ is any object of $\topos$, then $\rel(X)$ is an atomic modular quantale.
 }

\block{The notion of atomic quantale will be closely related to the notion of Hypergroupoid, which is a natural generalisation of the notion of canonical hypergroup which can be found in \cite{mittas1972hypergroupes} or \cite{Krasner_HG}.

\Def{A hypergroupoid $G$ is the data of:

\begin{itemize}
\item A set $E(G)$ of objects.
\item For each $e,e' \in E(G)$ a set $G(e,e')$ of ``arrows'' from $e$ to $e'$.
\item For each $g \in G(e_1,e_2)$ and $h \in G(e_2,e_3)$ an inhabited set $hg \subset G(e_1,e_3)$ of ``possible compositions''.
\item And for each $x \in G(e,e')$ an element $x^* \in G(e',e)$ called the inverse of $x$.
\end{itemize} 

With the following axioms:

\begin{itemize}

\item[(HG1)] $\forall e \in E(G), \exists 1_e \in G(e,e), \text{ such that } \forall x \in G(e',e), 1_e x = \{x\} \text{ and } \forall x \in G(e,e'), x 1_e = \{x\}. $

Such a $1_e$ is unique and is also denoted by $e$.

\item [(HG2)] for all $x,y,z$ three arrows such that the composition $xy$ and $yz$ are defined, one has $(xy)z = x (yz)$ where the product of an element $x$ with a set $S$ is defined by $xS = \bigvee_{s \in S } xs$.

\item [(HG3)] if $x \in yz$ then $z \in y^* x $ and $y \in x z^*$.

\end{itemize}

}

We will adopt the convention that if $g$ and $g'$ are two non-composable arrows of a hypergroupoid then $gg'$ is defined to be the empty set. Or constructively that for two general arrow $g$ and $g'$ of a hypergroupoid $gg' = \{u| g \text{ and } g' \text{ are composable and } u \in gg' \}$.

}

\block{\label{At_lem1} \Prop{Let $G$ be a hypergroupoid, and let $X$ be the set of all arrows:  
\[ X = \bigvee_{e,e' \in E(G)} G(e,e') \] 
Then $\powerob(X)$ endowed with the following structure:

\[ U^* = \{ x^*,x \in U \} \]

\[ UV = \left\lbrace t | \exists u \in U, v \in V, u,v \text{ are composable and } t \in uv \right\rbrace = \bigvee_{u \in U \atop v \in V} uv \]

is a (atomic) modular quantale.
}

\Dem{ $\powerob(X)$ is by definition an atomic locale. We check the remaining axioms:

\begin{itemize}

\item The associativity of the product: if $g \in U(VW)$ then $\exists u \in U$, $f\in (VW)$ such that $g \in uf$ but $b \in VW$ means that there exists $v \in V$, $w \in W$ such that $f \in v.w$. hence $g \in u(vw) = (uv)w$. So there exists an $h \in (uv)$ (in particular $h \in UV$) such that $g \in hw$ and, $g \in (UV)W$. The reverse inclusion is exactly the same.

\item The composition is clearly bi-linear because it is defined so that:

 \[ UV = \bigvee_{u \in U \atop v \in V} \{u\}.\{v\} \]

\item The Set $1= \{1_e, e \in E(G) \}$ is also clearly a unit for $Q$, because for any $u \in G(e,e')$, $\{u\}.\{1_e\} = \{u\}$ and $\{ u \} \{1_{e''} \} = \emptyset$ for any other $e''$,  hence $\{u\}1=\{u\}$, one obtains the general result by bi-linearity (and symmetry for the fact that $1$ is also a left unit).

\item $(UV)^* = V^* U^* $: Let $x \in V^* U^*$ i.e. $x \in v^* u^*$ for $u$ and $v$ respectively in $U$ and $V$. Then $v^* \in x .u $, $u \in x^* v^* $, and finally $x^* \in u v$, ie $x \in (UV)^*$. This reasoning can be conducted backwards to obtain the reverse inequality.

\item The modularity law:

Let $x \in U \wedge VW$ then $x \in U$ and $x \in vw$ with $v \in V$, $w\in W$. 
We have $v \in xw^*$ hence $v \in (UW^* \wedge V)$ and $x \in (UW^* \wedge V)W$.

\end{itemize}

}
}
\block{\label{At_lem2}\Prop{Conversely, any atomic modular quantale is of the form $\powerob(X)$ where $X$ is the set of arrows of a hypergroupoid. }

\Dem{ $Q$ be an atomic modular quantale. So $Q = \powerob(X)$ for some set $X$.

Let $E = \{x \in X| x \in 1_Q\}$. In order to simplify notations we will identify an element of $X$ with the corresponding singleton element in $Q$.

For any $q \in X$, as $1_Q q=q$ there exists $e \in E$ such that $e x =x$. Such an $e$ is unique because if  $e'x =x $ then $e' e x = x$. But as $e'$ and $e$ are subobjects of $1$ in a modular quantale, one has $ee'= {e}\wedge {e'}$ and in particular $x = (e' \wedge e)x$. Hence $e' \wedge e$ is inhabited and finally $e = e'$. Similarly, for each $x \in Q$ there is a unique $e' \in E$ such that $xe'=x$.

Let: 
\[ G(e,e') =\{ x \in X | xe=x \text{ and } e'x=x \} \]

We just show that $X$ is the disjoint union of all the $G(e,e')$. If $x \in G(e,e')$ then $x^* \in G(e',e)$. Also, if $a \in G(e,e')$ and $b \in G(e',e'')$ then $ba \subset G(e,e'')$, indeed, if $c \in ba$ then there exists a unit $f \in E$ such that $fc = c$. In particular $c \in fba = fe''ba=( f \wedge e'')ba$. Hence $e'' =f$.

We will prove that this is indeed a hypergroupoid structure:

\begin{itemize}

\item Let $a \in G(e,e')$ and $b \in G(e',e'')$. One has $e'a = a$ and $be'=b$.
hence $e' \in aa^*$ and $b \in baa^*$ so $ba$ is inhabited. 

\item if $a \in G(e',e)$ then $ea=a$ and $ae'=a$ by definition.

\item the associativity of the product comes from the associativity of the product of the quantale and the fact that if there exists $a \in uv$ for $u,v \in X$ then $u$ and $v$ are composable elements: this assert that the product of the hypergroupoid is exactly the product of the quantale, and its associativity follows by restriction to composable pairs of morphisms.

\item If we assume that $x \in yz$ then $x  = x \wedge yz \leqslant (xz^* \wedge y) z$
hence $y \wedge xz^*$ is inhabited so $y \in xz^*$. 

\end{itemize}

Finally as we have proved already that multiplication in $Q$ and in $X$ are essentially the same  it is a routine check to prove that $Q$ will be isomorphic to $\powerob(X)$ as a modular quantale.

}

}

\block{We have essentially proved that an atomic modular quantale is the same thing as a hypergroupoid. If we define a morphism of hypergroupoid from $G$ to $G'$ to be an application $f$ from $E(G)$ to $E(G')$ and a collection of maps (all called $f$) $f : G(e,e') \rightarrow G'(f(e),f(e'))$ such that $f(1_e)=1_{f(e)}$ and $f(xy) \subset f(x) f(y)$\footnote{Of course $f(xy)$ denote the direct image by $f$ of the set $xy$.},  then one has:

\Th{There is an anti-equivalence of category between the categories of atomic modular quantale (with weakly unital morphisms) and the categories of hypergroupoid.}

where ``weakly unital morphisms'' of (atomic) modular quantales are those defined in \ref{RelRep_morphism}, except that we no longer assume them to preserve the unit, but only to satisfy the weaker conditions $1 \leqslant f(1)$.

\Dem{Let $Q$ and $Q'$ be two atomic modular quantales and $g :Q' \rightarrow Q$ a morphism of modular quantales. By \ref{At_lem2}, $Q$ and $Q'$ can be written $Q'= \powerob(X')$ and $Q=\powerob(X)$ where $X'$ and $X$ are the sets of all arrows of two hypergroupoids $G$ and $G'$. As $g$ is in particular a morphism of locale, it induces an application $f : X \rightarrow X'$ characterised by the fact that for all $x\in X$, $f(x)$ is the unique element of $X'$ such that $x \in g(f(x))$.

In particular, let $c \in ab$ then $ab \subset g(f(a))g(f(b))$ hence $c \in g(f(a)f(b))$ ie $f(c) \in f(a)f(b)$. This proves that $f(ab) \subset f(a)f(b)$. As $1 \leqslant g(1)$, and $1\in G$ corresponds to $E(G) \subset X$, the application $f$ acts on the unit set and preserves the identity element. One also has $x^* \in g(f(x))^*=g(f(x)^*)$ hence $f(x)^*=f(x^*)$ and finally, if $g \in G(e,e')$ then $e\in g^*g$, $e'\in gg^*$ and $f(e)\in f(g)^*f(g)$, $f(e') \in f(g)f(g)^*$ which proves that $f(g)$ is an element of $G(f(e),f(e'))$ which concludes the proof that $f$ is a morphism of hypergroupoid.

Conversely, if $f$ is a morphism of hypergroupoids, then as $f(1_e) \in f(g)f(g^*)$ one can conclude that $f(g^*) \in f(g)^*f(1_e) = f(g)^*$ hence $f(g^*)=f(g)^*$. One can then define $g=f^{-1}$ which is a frame homomorphism and compatible with multiplication and involution, and $1 \leqslant g(1)$ because each unit is sent by $f$ to a unit, hence it is a morphism of modular quantales. These two constructions are clearly compatible with compositions and inverse from each other, hence, together with propositions \ref{At_lem1} and \ref{At_lem2}, this concludes the proof of the theorem.

}

}

\block{ Finally we investigate the case of Grothendieck quantale: 

\Prop{Let $Q$ be a modular quantale, and $G$ be the corresponding hypergroup. The following conditions are equivalent:

\begin{enumerate}

\item $Q$ is a Grothendieck quantale

\item for every arrow $f$ in $G$ there exists two arrows $u$ and $v$ such that $f=uv^*$ with $uu^*$ and $vv^*$ units.

\end{enumerate}
}

\Dem{If $Q$ is a Grothendieck quantale, then (by $(Q10)$) any element of $Q$ can be written as a supremum of elements of the form $uv^*$ with $uu^*\leqslant 1$ and $vv^* \leqslant 1$, as $Q$ is atomic we can write these $u$ and $v$ as union of atoms and hence any elements of $Q$ can be written as a supremum of element of the form $uv^*$ where $uu^*$ and $vv^*$ are units. In particular, any $f \in G$ is an atom of $Q$, and hence should be of the form $uv^*$.

\bigskip

Conversely, if $G$ satisfies condition $(2)$ then any element of $Q$ can be written as a union of its atoms, which are all of the form $uv^*$ with $uu^*$ and $vv^*$ units (and hence $\leqslant 1$).

}

\Def{An element $f$ of a hypergroupoid such that $ff^*$ is a unit is called a simple element. An element which can be written in the form $fg^*$ with $f$ and $g$ simple is called a semi-simple element. A hypergroupoid satisfying the condition of the proposition, i.e. such that every element is semi-simple, will be called a semi-simple hypergroupoid. }

Hence an atomic Grothendieck quantale is essentially the same thing as a semi-simple hypergroupoid.

}

\subsection{Hypergroupoid algebra}
\label{hga}

\blockn{In this section we consider an atomic topos $\topos$, an arbitrary object $X$, the quantale $Q$ of relations on $X$ and the corresponding hypergroupoid $G$. We assume that $G$ (the set of all arrow of $G$) is decidable. This implies that $X$ is a decidable object, indeed:

\Lem{Let $G$ be a decidable hypergroupoid, then its set of units is complemented. And if $G$ is associated to an object $X$ of an atomic topos then $X$ is decidable.}
\Dem{Let 
\[\Delta^c = \{g \in G | \forall e \in E(G), e \neq x \}. \]
We will prove that $\Delta^c$ is a complement of $E(G)$.
They are disjoint and for every $g \in G$ there exists a unique $e$ such that $eg=g$. As $G$ is decidable, either $g=e$ or $g \neq e$. If $g=e$ then $g \in E(G)$. If $g \neq e$ then for all $e'$, one has $e'=g \Rightarrow e'=e$ because of the uniqueness of $e$, and hence $e' = g$ yields a contradiction, so $g \in \Delta^c$.

In particular, as $\powerob(G)$ is isomorphic to $\sub(X \times X)$, and as the diagonal sub-object of $X$ corresponds to the set of units of $G$, this proves that $X$ is decidable.

}

In this situation, the convolution product defined in section \ref{ConvProd} gives a convolution on functions on $G$ with value in $\Rlscp$. We will give necessary and sufficient conditions in order that the convolution product induces an interesting multiplication on some algebra.

As, in this situation, the convolution product depends on $\topos$ and not only on $G$, these conditions will be expressed in terms of the logic of $\topos$. In the next section we will focus on the case of a semi-simple hypergroupoid, in this situation the topos $\topos$ will be  canonically determined and it will be possible to reformulate the definition given here in more explicit terms. 
}

\block{\label{hga_card}We will need a few generalities about cardinals of sets in a constructive setting in order to be able to give an internally valid proof of the mains results of this section.

\Def{Let $X$ be a decidable set, then the cardinal of $X$ is defined by:

 \[ |X| = \left( \sum_{x \in X} 1 \right) \in \Rlscp \] }

We remind the reader that $\Rlscp$ contains an element $+\infty$.
The following lemma gives two properties that completely characterise the cardinal of a set.

\Lem{ \begin{itemize} 

\item For $n \in \mathbb{N}$ one has $n \leqslant |X|$ if and only if there exists $x_1, \cdots x_n$ pairwise disjoint elements of $X$. 

\item For $q \in \mathbb{Q}$, $q < |X|$ if and only if there exists an $n \in \mathbb{N}$ such that $q<n \leqslant |X|$.

\end{itemize} } 
 
\Dem{ Let $q \in \mathbb{Q}$, such that $q < |X|$. By definition, there exists $ x_1,\dots,x_n \in X$ pairwise distinct such that $\exists q_1, \dots, q_n<1$ with $q< \sum q_i$. This can be rewritten as $\exists x_1,\dots,x_n \in X$ such that $q<n$. This proves the second point of the lemma assuming the first.

If there are $n$ distinct elements in $X$, any $q <n$ is also smaller than $|X|$ hence $n \leqslant |X|$. Conversely, if $n \leqslant |X|$ then $(n-\frac{1}{2}) < |X|$, so there is an integer $m$ with $(n-\frac{1}{2})<m$ and $x_1, \dots x_m$ distinct element in $X$. As $n \leqslant m$ one also has $n$ distinct elements, this concludes the proof of the first point of the lemma.

}
 
\Prop{If $X$ is a decidable set, the following conditions are equivalent:

\begin{enumerate}
\item $X$ is finite.

\item $|X|$ is an integer.

\item $|X|$ is a (finite) continuous real number (ie an element of $\Rt$).
\end{enumerate}

}
For an example of a set $X$ with $|X|<\infty$ but not satisfying these properties, one can take any non-complemented sub-set of a finite (decidable) set.
\Dem{$1. \Rightarrow 2.$ is clear because a finite decidable set is isomorphic to $\{1, \dots, n \}$ and hence as cardinal $n$.

$2. \Rightarrow 3.$ is also clear.

Assume $3.$, then there exist $q,q'$ such that $|q-q'|<\frac{1}{2}$ and $q<|X|<q'$. There exists an integer $n$ such that $q <n \leqslant |X| < q'$, and $x_1, \dots, x_n$ pairwise distinct elements of $X$. Let $x \in X$ then there are two possible cases: either $x=x_i$ for some $i$, or $x$ is distinct from all the other $x_i$ (this is proved recursively on $n$ using the decidability of $X$). But if $x$ is distinct from all the $x_i$ then $(n+1) \leqslant |X|$ and $q \leqslant n < (n+1) \leqslant q$ which yields a contradiction. So $x = x_i$ for some $i$, and $X$ is indeed finite.

}

We also note that the same argument yields the following result: If $X$ is decidable, and we have a function p:$X \rightarrow \mathbb{N}^{>0}$ such that $\sum_{x \in X} p(x)$ is an integer, then $X$ is finite.

}

\block{\label{hga_hgacoef} Let $g$ and $g'$ be two arrows in $G$, and $[g]$,$[g']$ the characteristic functions of the singletons $\{g\}$ and $\{g'\}$ then:

\[ ([g]*[g'])(x,y) = \sum_{  z  } [xgz \wedge z g'y] =  |\{z| x g z \text{ and } z g' y \}| \]

\Def{ Let $\hgacoef{a}{g}{g'}$ denote the evaluation in $a$ of the function $[g]*[g']$.

We also define for $g \in G(e,e')$, 

\[ \begin{array}{c c}
|g|_l = \hgacoef{e}{g^*}{g} & |g|_r = \hgacoef{e'}{g}{g^*}
\end{array} \]
} 

\Prop{ $\hgacoef{a}{g}{g'}$, $|g|_l$ and $|g|_r$ can be computed internally using the formulas: 

For any $u \in e$ the source unit of $g$:
\[ |g|_l = |\{z | z g u\}|. \]
For any $v \in e'$ the target unit of $g$:
\[ |g|_r = |\{z| v g z \}| \]
For any $(x,y) \in a$:
\[ \hgacoef{a}{g}{g'} = |\{ t | x g t \text{ and } t g' y \}| \]
}

\Dem{The formula for $\hgacoef{a}{g}{g'}$ is essentially its definition. The two other formulas follow easily.
The fact that each time the value internally does not depend on any choice of (internal) elements is clear because the various possible choices all belong to a same atom. }

One also mentions the two easy (but important) relations: 

\[ \begin{array}{c c} |g^*|_l = |g|_r,   & \hgacoef{a}{g}{g'}= \hgacoef{a^*}{g'^*}{g^*} \end{array}\]

}

\block{\label{hga_murel}
Here are some of the important combinatorial properties of these coefficients:

\Th{ For all pair $g,g'$ of composable arrows and all $a \in gg'$ on has:
\begin{enumerate}
\item $ \displaystyle \hgacoef{a}{g}{g'}|a|_l = \hgacoef{g'}{g^*}{a}|g'|_l$
\item $ \displaystyle \hgacoef{a}{g}{g'}|a|_r = \hgacoef{g}{a}{g'^*}|g|_r$
\item $\displaystyle |g|_l |g'|_l = \sum_{a \in gg'} \hgacoef{a}{g}{g'} |a|_l $
\end{enumerate}
}

\Dem{ let $e_1$,$e_2$ and $e_3$ be the units such that $g' \in G(e_1,e_2)$ and $g \in G(e_2,e_3)$.

Let $e \in e_1$ be arbitrary then let:

\[ X_a = \{(u,v) | u g v \text{ and } v g' e \text{ and } u a e \} \]
 \[ X = \{ (u,v) | u g v \text{ and } v g' e \} = \bigsqcup_{a \in gg'} X_a. \]

The cardinality of $X_a$ can be computed in two different ways, on one side:

\[ |X_a| = \sum_{u \text{ s.t. } u a e} |\{v | u g v \text{ and } v g' e \} = \sum_{u \text{ s.t. } u a e} \hgacoef{a}{g}{g'}=|a|_l \hgacoef{a}{g}{g'}. \]

On the other side:

\[ |X_a| = \sum_{v \text{ s.t. } v g' e } |\{u| v g^* u \text{ and } u a e \}| = \sum_{v \text{ s.t. } v g' e } \hgacoef{g'}{g^*}{a} =|g'|_l \hgacoef{g'}{g^*}{a}. \]

The equality of the two results gives $(1.)$. the result $(2.)$ is the dual (one can use a similar proof or apply $*$ every where).

Similarly,

\[ |X| = \sum_{v \text{ s.t. } v g' e} |\{u|u g v \}| = \sum_{v, v g' e}  |g|_l = |g|_l|g|_r \]

Hence $3.$ comes from the fact that $X$ is the disjoint union of the $X_a$.

}
} 

\block{
\Th{Let $G$ be a decidable hypergroupoid represented in a topos i.e. corresponding to the modular quantale $\rel(X)$ for some object $X$ of a topos. The following propositions are equivalent: 
\begin{itemize}
\item for all $g \in G, g:e' \rightarrow e$ the value of $[g]*[g^*]$ at $e$ is a finite continuous real number.
\item for all $g,g' \in G$ the set $gg'$ is a finite set and one has a formula of the form:

\[ [g]*[g'] = \sum_{a \in gg'} \hgacoef{a}{g}{g'} [a] \]

with $\hgacoef{a}{g}{g'}$ (finite) positive integers.
\end{itemize}
In the case where this conditions are verified we say that $(G,X)$ is locally finite.
}

\Dem{The second condition clearly implies the first. We assume the first condition.

This means that for all $g \in G$ one has internally $\forall y \in X$ the set $\{x | x g y \}$ is a finite  set (its cardinal is $\hgacoef{e}{g}{g^*}$ and is continuous by assumption). In particular $\{x | x g y \}$ as a finite subset of a decidable set is complemented in $X$, hence $\{z| x g z \text{ and } z g' y \}$ is finite as a complemented subset of a finite set.

This proves that the evaluation of $[g]*[g']$ at every point is indeed a positive integer. All that remains to do is check  that $gg'$ is finite but this comes from the fact that all coefficients appearing in the sum of lemma \ref{hga_murel} ($3.$) are strictly positive integers, hence the sum has to be indexed by a finite set (see the remarks at the end of \ref{hga_card}).

}

The situation described in the second condition is basically the best we can hope: we get a $\mathbb{Z}$ algebra generated (as a group) by the symbol $[g]$ for $g \in G$. We will call this algebra $A_G$.

But conversely, if we want to have any interesting convolution structure coming from the construction done in \ref{ConvProd} on a set of functions on $G$ with value in continuous numbers, we need to have the first condition. This proves that the ``locally finite'' hypothesis is exactly the good hypothesis for getting an interesting convolution product.

}

\subsection{Semi-simple Hypergroupoid algebra}
\label{sshga}

\blockn{In this section we assume that $X$ is a bound of $\topos$. The hypergroupoid $G$ is now semi-simple, and $\topos$ is fully determined by $G$, so we should be able to express the value of $\hgacoef{a}{g}{g'}$ in terms of the structure of $G$. }

\block{\label{sshga_main}\Prop{In a semi-simple decidable hypergroupoid one has:
\[ \hgacoef{a}{g}{g'} = \sup_{a=xy^* \atop x,y \text{ simples}} |g^*x \wedge g'y | \]
}

The $\sup$ in the proposition is taken in $\Rlscp$, this means in particular that the coefficient is an integer if and only if the supremum is reached.

\Dem{ Let $q < \hgacoef{a}{g}{g'}$ i.e. $q<n \leqslant \hgacoef{a}{g}{g'}$.
This means that for every $(x,y) \in a$, there exists $(v_1, \dots v_n)$ pairwise distinct in $B$ such that for all $i$, $(x g v_i)$ and $(v_i g' y)$.

This means that there is a surjection $(x,y): t \twoheadrightarrow a$ and a collection of $n$ maps $v_1, \dots v_n$ from $t$ to $B$ pairwise distinct\footnote{As we will soon assume that $t$ is an atom the precise meaning of `distinct' is not important.}, such that for all $i$, $(x,v_i)$ has value in $g$ and $(v_i,y)$ has value in $g'$.

If we choose any atom on $B$ which maps to $t$, the composite is still a surjection on $a$. Hence we can freely assume that $t$ is a unit of $G$, and that $x,y,v_1, \dots v_n$ are arrows in $G$. And the two relations $x(t) g v_i(t)$ and $v_i(t) g' y(t)$ become: $v_1, \dots v_n \in g'y \wedge g^* x$ and $x(t) a y(t)$ became $a=xy^*$, so 
\[ q< n \leqslant \sup_{a=xy^*} |g'y \wedge g^*x|. \]

Conversely, if $q< \sup_{a=xy^*} |g'y \wedge g^*x|$, then for some $x,y$ simple such that $a=xy^*$,

\[ q < n \leqslant |g'y \wedge g^*x|. \]

So there exist $v_1,\dots,v_n$ pairwise distinct in $g'y \wedge g^*x$. If $n>0$, this implies that $g,g'$ is composable (if $n=0$, then $q$ is smaller than any cardinal). Let $e$ be the target of $g$ and the source of $g'$.

In this situation, for any $u \in e$,  $v_1(u) \dots v_n(u)$ are $n$ pairwise distinct elements in $\{z| x(u) g z \text{ and } z g' y (u) \}$ and hence $q<n \leqslant \hgacoef{a}{g}{g'}$.

And this concludes the proof.
}
}

\block{\Prop{An element $g \in G$ is simple if and only if $|g|_l=1$.}
\Dem{Let $g \in G(e,e')$. 
Internally one has by the proposition \ref{hga_hgacoef}:
 \[ \forall x \in e, |g|_l=| \{z |z g x \}| \] 
Hence $|g|_l=1$ exactly means that for all $x$ there exists a unique $z$ such that $z g x$ which means that $g$ is a partial function, i.e. a simple element.
}
}

\block{\Prop{Let $G$ be a semi-simple decidable hypergroupoid, then $g \in G$ is left finite if and only if there exists a simple arrow $u$ such that $(g,u)$ is composable , $gu$ is finite and contains only simple elements. Moreover in this case $|g|_l$ is the cardinal of the set $(gu)$.}

\Dem{Translating into an external result the internal formulation of the left finiteness of $g$ would actually give us exactly the statement of this theorem, but this translation require some work (essentially done in the proof \ref{sshga_main}) which can be avoided by the use of the combinatorial identities we already proved.

Assume first that $gu$ is finite and contains only simple elements $\{x_1,\dots,x_n \}$. then by the formula $3.$ of \ref{hga_murel} and replacing by $1$ the left cardinal of simple elements, one gets that:

\[ |g|_l = \sum_{i=1}^n \hgacoef{x_i}{g}{u} \]

But one can see on the formula given in the proposition \ref{hga_hgacoef} that $\hgacoef{x_i}{g}{u} \leqslant |u|_l =1$, hence $\hgacoef{x_i}{g}{u} =1$ and $|g|_l=n$ which implies that $g$ is left finite.

Conversely, assume that $|g|_l=n$ for some $n$. By \ref{sshga_main} one has:

\[ |g|_l = \sup_{u \text{ simple}}|gu| \]

In particular as $(n-1/2) < |g|_l$, there exists $u$ simple such that $|gu|=n$. one has then:

\[ |g|_l|u|_l = \sum_{x \in gu} \hgacoef{x}{g}{u} |x|_l \]
\[ n= \sum_{i=1}^n  \hgacoef{x_i}{g}{u} |x_i|_l \]

This implies first that all $\hgacoef{x_i}{g}{u} |x_i|_l$ have an opposite, hence they are all continuous numbers, and hence integers. Moreover as all the  $\hgacoef{x_i}{g}{u} |x_i|_l$ are $\geqslant 1$ they have to be all equal to $1$, hence all the $x_i$ are simple and this concludes the proof.

}

}

\subsection{The category of $\topos$-groups}
\label{Tgroup}
\blockn{In this section, we will show that in the locally finite case, the algebra $A_G$ we obtained can be seen as a particular subalgebra of endomorphisms of the free group $\mathbb{Z}B$ generated by $B$ in the logic of $\topos$. This gives an abstract interpretation of the algebra $A_G$.

We will also show that in the semi-simple case the category of groups of $\topos$ embeds as a full subcategory of the category of $A_G$-modules, and that this embedding induces an equivalence between $\Qt$-vector spaces in  $\topos$ and full $A_G \otimes \mathbb{Q}$-modules.
 }

\block{\label{Tgroup_endoZB}Let $E=\mathbb{Z}B$ be the free group generated by $B$ in $\topos$ ($E$ is a group object of $\topos$). In particular $E \otimes E = \mathbb{Z}(B \times B)$ is the free group generated by $B \times B$, hence as $B$ is assumed to be decidable one can define a bi-linear map $\Delta: E \times E \rightarrow \Zt$ which sends $(b,b')$ to $1$ if $b=b'$ and to $0$ if $b \neq b'$. 

Let $f$ be an endomorphism of $E$. One can associate to $E$ the function of ``matrix elements'' of $f$,  $\rho(f): (b,b') \mapsto \Delta(b,f(b'))$ from $B \times B$ to $\mathbb{Z}$. The map $\rho: f \mapsto \rho(f)$ is injective, and one has:

\[f(b) = \sum_{b'} \rho(f)(b',b) b' \]

Also, $\rho(f \circ f') = \rho(f) * \rho(f')$ for the convolution of functions on $B \times B$ (here the sum involved in the convolution product will be finite, hence it is defined for integer valued functions). So we just have to understand the image of $\rho$:

\Prop{A function $f$ from $G$ to $\mathbb{Z}$ belongs to the image of $\rho$ if and only if it verifies the following two properties:

\begin{enumerate}

\item If $f(g) \neq 0$ then $g$ is left finite.

\item For each unit $e$ of $G$, there is at most a finite number of arrows $g \in G$ pointing to $e$ such that $f(g)$ is non zero.

\end{enumerate}

 }

In particular, property $(1.)$ tells us that if we are not in the locally finite case, then the algebra of group homomorphisms of $E$ is in some sense `too small'.

Also, in the locally finite case, the algebra $A_Q$ is identified with the sub-algebra of endomorphisms of $E$ such that $f(g)$ is non zero only for a finite number of $g \in G$.

\Dem{Internally, a function $f: X \times X \rightarrow \Zt$ corresponds to a group homomorphism if and only if for all $x \in X$ there is only a finite set of $x'$ such that $f(x',x)$ is non zero. Indeed, the corresponding group homomorphism has to send $x$ to $\sum_{x'} f(x',x).x'$.

The cardinality of the set of $x'$ such that $f(x',x)$ is non zero defines a function $\tau$ on $X$, whose value at any atom $e$ of $X$ (ie at any unit of $G$) is given by:

\[ \tau(e) = \sum_{g \in G(e',e), f(g) \neq 0} |g|_l \]

Indeed, for any $x \in X$ the set of $x'$ such that $f(x',x)$ is non zero is partitioned by the various $g \in G(e',e)$ such that $(x',x) \in g$ and each of this set has cardinality $|g|_l$ (because the value of $f(x',x)$ only depends on the atom that contains it).

So $\tau(e)$ is an integer if and only if each of the $|g|_l$ are integers (this is condition $1.$) and if they arise in finite number (this is condition $2.$).

}

}

\blockn{In the rest of this section (and also in the rest of this paper), we assume that the representation of $G$ in $\topos$ is locally finite.}

\block{\label{Tgroup_fct}Let $F$ be any group object of $\topos$. Then $\hom(E,F)$ is a right $\hom(E,E)$ module (where the $\hom$ denotes the internal group homomorphisms). The units of $G$ act as family of disjoint projections on $E$, and hence on $\hom(E,F)$. Let:

 \[ \tilde{F} = \bigoplus_{e \in E(G)} Hom(E,F).e \]

equivalently, $\tilde{F}$ is the subset of $\hom(E,F)$ of elements $x$ such that there exists a finite set $I$ of units such that $x.e_I =x$ where $e_I = \sum_{e\in I} e$.

\Prop{ For all $F$, $\tilde{F}$ is a full\footnote{A module $M$ over a non unitary ring $A$ is said to be full if the map $A \times M \rightarrow M$ is surjective.} right $A_Q$-module.

This gives a functor from group objects of $\topos$ to full right $A_Q$-modules.

Also $\tilde{E}$ is $A_Q$ seen as a right $A_Q$-module.}

\Dem{Clear from the observation that $\tilde{F}$ is the subset of $\hom(E,F)$ of elements $x$ such that there exists a finite set $I$ of units such that $x.e_I =x$ where $e_I = \sum_{e\in I} e \in A_Q$. And $\tilde{E}$ identify with $A_Q$ thanks to \ref{Tgroup_endoZB}}

Actually, $\tilde{F} = \hom(E,F).A_Q$.

}

\block{\label{Tgroup_equiv}Assume that $G$ is semi-simple, so $B$ is a bound of $\topos$ and the category of atoms of $B$ and morphisms between them (i.e. the category of units of $G$ and simple arrows between them) endowed with the atomic topology is a site of definition of $\topos$. 

If $F$ is a group object of $\topos$ then for each $e$ atom of $B$, $F(e)=\tilde{F}.e$ and the action of a simple arrow $f$ from $e$ to $e'$ is given be the action of $[f]$ on $\tilde{F}$. In particular, the sheaf corresponding to $F$ is fully determined by $\tilde{F}$ and any $A_Q$-linear morphism from $\tilde{F}$ to $\tilde{F'}$ gives rise to a morphism of sheaves and one can conclude that:

\Lem{When $G$ is semi-simple (and locally finite) the functor from $\topos$-groups to $A_G$-modules defined in \ref{Tgroup_fct} is fully faithful.}

Unfortunately, if we start from a general $A_G$-module we only get a pre-sheaf over the site of units. In the general case, we have not found a characterisation of the $A_G$-modules corresponding to $\topos$-group simpler than the definition of a sheaf. But in the case where we assume that all the coefficients $|g|_l$ are invertible, then the action of the $[g]$ for non simple $g$ will automatically  turn our pre-sheaf into a sheaf.

\Prop{(Still under the assumption that $G$ is semi-simple and locally finite) Let $M$ be an $A_G$ module such that for every $g \in G$ the integer $|g|_l$ acts (by multiplication) as a bijection on $M$. Then $M$ comes from a $\topos$-group.

In particular there is an equivalence of categories between $\Qt$-vector spaces and full right $A_G \otimes \mathbb{Q}$-modules.
}

\Dem{We will check that under this assumption the pre-sheaf of $\tilde{M}.e$ is actually a sheaf for the atomic topology.

Let $e$ be a unit of $G$, let $f$ be any simple arrow starting at $e$. And let $m \in M.e$ such that for any to simple arrows $g,h$ targeting $e$ such that $fg=fh$ one has $mg=mh$. We need to prove that there exists a unique $n$ such that $nf=m$.

The uniqueness is easy: if $m=n.f$ then $m.f^* = n.ff^* = |f|_r n$ so as $|f|_r$ is invertible on has $n = \frac{1}{|f|_r} m.f^*$.

Conversely, we will prove that $n = \frac{1}{|f|_r} m.f^*$ provides a solution.

\[ n.f = \frac{1}{|f|_r} m.f^*f =\frac{1}{|f|_r} \sum_{a \in f^*f} \hgacoef{a}{f^*}{f} m.a \]

let $a \in f^*f$. $a$ can be written in the form $a=gh^*$ with $g$ and $h$ simple.

\[ gh^* \in f^*f \Rightarrow f \in fgh^* \Rightarrow fg \in fh \]

Hence one has $fg=fh$ and by the assumption on $m$, $m.h=m.g$.

the relation $ m.g = m.h$ implies,

\[ m.[g][h^*] = m.[h][h^*] \]

\[ [h][h^*]=|h|_r.[e]\]

By \ref{hga_murel} one has:

\[ [g][h^*] = \frac{|h|_r}{|a|_r} [a]\]

so $m.[g][h^*] = m.[h][h^*]$ becomes: 
\[m.a = |a|_r .m \]

and we can conclude that: 

\[n.f = \frac{1}{|f|_r} \sum_{a \in f^*f} \hgacoef{a}{f^*}{f} |a|_r m = |f^*|_r m = m\]

again by \ref{hga_murel} and the fact that $f$ is simple hence $|f^*|_r = |f|_l=1$.

}

}

\subsection{$A_G$ as a quantum dynamical system}
\label{Qsyst}
\blockn{In this section we construct the regular representation of $A_G$. We show that the $C^*$ algebra generated by $A_G$ comes with a canonical action of $\mathbb{R}$. There is also a regular representation of $A_G$, attached to a $KMS_1$ state and defining a $C^*$ algebra $\Cred(G)$ by completion.}

\block{\label{Qsyst_adj}Let $\mathcal{H}$ be a real Hilbert space in $\topos$. Then

 \[ \tilde{\mathcal{H}} = \bigoplus_{e \in E(G)} \hom(e,\mathcal{H}) \]

is an $A_G$ vector space. The internal scalar product on $\mathcal{H}$ gives rise to a scalar product on each of the $\hom(e,\mathcal{H})$ and turns $\tilde{\mathcal{H}}$ into a pre-Hilbert space.

\Prop{In the action of $A_G$ on $\tilde{\mathcal{H}}$ one has:

\[ [g]^* = \frac{|g|_l}{|g|_r}[g^*] \] 

And $[g]$ has norm smaller than $|g|_l$.

}

\Dem{ Let $g \in G$, $v \in \hom(e,\mathcal{H})$ and $v' \in \hom(e',\mathcal{H})$.

If $g$ is not an arrow from $e$ to $e'$, then both $ \scal{v}{v'[g]}$ and $ \scal{v[g^*]}{v'}$  are zero (hence equal). If $g$ is an arrow from $e$ to $e'$, then for any $x \in e$ one has:

\[ \scal{v}{v'[g]}(x) = \scal{v(x)}{\sum_{ y, y g x} v'(y) } =\sum_{y,y g x} \scal{v(x)}{v'(y)}\]

But $g$ is an atom of $B \times B$, and $\scal{v(x)}{v'(y)}$ is a function on $B \times B$ so its value does not depend on $x$,$y$ as long as they belong to $g$. Hence, for any $(x,y) \in g$ one has:

\[ \scal{v}{v'[g]} = |g|_l \scal{v(x)}{v'(y)} \]

Similarly:

\[ \scal{v[g^*]}{v'} = |g|_r \scal{v(x)}{v'(y)} \]

Finally: 
\[ \scal{v}{v'[g]}= \frac{|g|_l}{|g|_r} \scal{v[g^*]}{v'} \]

And the first result follows.

The second result follows from 
\[\scal{v}{v'[g]} = |g|_l.\scal{v(x)}{v'(y)}\leqslant |g|_l \Vert v \Vert \Vert v' \Vert .\]

As $\Vert v(x) \Vert = \Vert v \Vert$.

}
}

\block{\label{Qsyst_chirel}\Def{For $g \in G$, We will denote:

\[ \chi(g) = \frac{|g|_l}{|g|_r} \] 

}
One has:

\[ \chi(g^*) = \chi(g)^{-1} \]
\[ [g]^* = \chi(g) [g^*] \]

and also the more surprising result: 

\Prop{For any three arrows $a,g,g'$ of $G$ such that $a \in gg'$:

\[ \chi(a)=\chi(g)\chi(g') \] }

\Dem{We will just need several applications of the first two points of theorem \ref{hga_murel}.

\[ \chi(a) = \frac{|a|_l}{|a|_r} = \frac{|a|_l \hgacoef{a}{g}{g'}}{|a|_r \hgacoef{a}{g}{g'} } = \frac{ |g'|_l }{|g|_r} \frac{\hgacoef{g'}{g^*}{a} }{\hgacoef{g}{a}{g'^*}}  \]

but:

\[\frac{|g'|_r}{|g|_l} \frac{\hgacoef{g'}{g^*}{a} }{\hgacoef{g}{a}{g'^*}} = \frac{|g|_l \hgacoef{g^*}{g'}{a^*}}{|g'|_r \hgacoef{g'^*}{a^*}{g} } = \frac{|g|_l\hgacoef{g}{a}{g'^*}}{|g'^*|_l  \hgacoef{g'^*}{a^*}{g}} = 1 \]

So we can conclude that:
\[ \chi(a) = \frac{ |g'|_l }{|g|_r} \frac{\hgacoef{g'}{g^*}{a} }{\hgacoef{g}{a}{g'^*}} =\frac{ |g'|_l }{|g|_r}\frac{|g|_l}{|g'|_r} = \chi(g) \chi(g') \] 
}
}

\block{\Def{Let $A_{G,\mathbb{Q}}$ be the algebra $A_G \otimes \mathbb{Q}$. It is endowed with the involution $\_^*$ defined by \ref{Qsyst_adj}, i.e.:

\[ [g] ^* = \chi(g) [g^*] \]

We also define the elements:

\[ e_g = \frac{1}{|g|_l} [g] \in A_{G,\mathbb{Q}} \]

which are additive generator such that $(e_g)^* = e_{g^*}$.

}

\Prop{If $G$ is semi-simple and locally finite, the functor which sends a $\topos$-Hilbert space $\mathcal{H}$ to the completion of $\tilde{H}$ is (one half of) an equivalence of categories between the category of internal Hilbert spaces of $\topos$, and the full\footnote{Here full, mean that $\mathcal{H}.A_{G,\mathbb{Q}}$ is dense in $\mathcal{H}$. } right Hilbert $*$-representations of $A_{G,\mathbb{Q}}$. }

\Dem{The proof is really similar to the case of $\Qt$-vector spaces done in \ref{Tgroup_fct} and \ref{Tgroup_equiv}. If we start from a $\topos$-vector space, we already proved that $A_{G,\mathbb{Q}}$ acts on $\tilde{H}$ by bounded morphisms and in a way compatible with the involution. So it extends to a full Hilbert $*$ representation of $A_{G,\mathbb{Q}}$ on the completion of $H$.

\bigskip

In the other direction, the sheaf of complex numbers on the site of units of $G$ is the constant sheaf. Hence if $\mathcal{H}$ is a Hilbert $*$ representation of $A_{G,\mathbb{Q}}$ then  $e \rightarrow \mathcal{H}[e]$ defines a pre-sheaf of $\Ct$-modules that we will denote by $H$.
The pre-sheaf $H$ is a sheaf by \ref{Tgroup_equiv}.

For every simple arrow $f:e \rightarrow e'$ one has:

 \[ [f].[f]^* = [f][f^*]\frac{1}{|f|_r} = [e']\frac{|f|_r}{|f|_r}=[e'] \]

and the induced map $\mathcal{H}[e'] \rightarrow \mathcal{H}[e]$ (i.e. the structural map of $H$) is an isometric injection. This proves that the scalar product $\mathcal{H}[e] \times \mathcal{H}[e] \rightarrow \mathbb{C}$ is in fact a morphism of sheaves $H \times H \rightarrow \Ct$ and hence this endows $H$ with an internal scalar product.
 
It remains to show that $H$ is internally complete. Let $\tilde{H}$ be its completion, let $h \in \tilde{H}(e)$.
 
Then, by (internal) density of $H$ in $\tilde{H}$, for every $ n \in \mathbb{N}$, there exists $f: e' \rightarrow e$, and $h' \in H(e')$ such that $\Vert h' - h.f \Vert<1/n$. But one can write:

\[\Vert h'[f]^* - h \Vert <1/n  \]

and $h'[f]^* \in H(e)$, hence $h$ can be approximated by elements of $H(e)$. As $H(e)$ is complete, this proves that $h \in H(e)$ and hence that $H$ is internally complete.

Finally this is an equivalence, because if we start from a full Hilbert $*$-representation $H$ of $A_{G,\mathbb{Q}}$ then the construction we just made corresponds to that of \ref{Tgroup_equiv} applied to $H.A_{G,\mathbb{Q}}$ hence as we applied a completion at the end, we will get $H$ back because $H.A_{G,\mathbb{Q}}$ is dense in $H$ by assumption.

}

}

\block{At this point we can either consider the closure of $A_{G,\mathbb{C}} = A_G \otimes \mathbb{C}$ in a specific representation: the one corresponding to the internal Hilbert space $l^2(B)$ of square summable functions on $B$, defining a $C^*$ algebra $\Cred(G)$, or we can take the universal $C^*$ algebra generated by $A_{G,\mathbb{Q}}$ that we will denote by $\Cmax(G)$. 

Both these algebras come with a time evolution $(\sigma_t)_{t \in \mathbb{R}}$ given by:

\[\sigma_t([g]) = \chi(g)^{it}[g] \]

This is a morphism of algebras because of \ref{Qsyst_chirel}.

Let $e$ be an atom of $B$. Then one has a map $e \hookrightarrow B$ which gives rise to a map $e \rightarrow l^2(B)$ and hence to a vector of the corresponding representation of $\Cred(G)$ that we will simply denote $l^2$.
An easy computation shows that the state on $\Cred$ induced by this vector is: 
\[ 
\eta_e([q])= \left\lbrace \begin{array}{l} 
1 \text{ if $q=e$} \\
0 \text{ otherwise.} 
\end{array} \right.
\]

If the set of units of $G$ is finite then the (renormalized) sum of all the $\eta_e$ is a state, in general we can define it without renormalization as a semi-finite weight (it is finite on the algebra $A_G$), we denote it by $\eta$.

\Prop{The GNS representation induced by $\eta$ is the $l^2$ representation, and $\eta$ verifies the KMS condition at temperature one. }

\Dem{The first part is clear: the GNS representation induced by $\eta$ is included in $l^2$ and contains all the vectors corresponding to the $e \in G$ (indeed, $[e]$ gives rise to this vector through the GNS construction). If it were a strict sub-representation then it would correspond internally to a sub Hilbert space of $l^2(B)$ containing all the basis vectors, which is impossible. 

For the second part:

\[ \eta([q] \sigma_{i}([q'])) = \chi(q')^{-1} \eta([q][q']) \]

If $q' \neq q*$ then both $\eta([q][q'])$ and $\eta([q'][q])$ are zero (because $e \in qq' \Rightarrow q' = q^*$). If $q'=q^*$ then

\[ \eta([q] \sigma_{i}([q']) = \chi(q) \eta([q][q^*])= \chi(q)  \hgacoef{e}{q}{q^*} = \chi(q) |q|_r = |q|_l = \eta([q'][q]) \]

}
}

\subsection{The time evolution of an atomic locally separated topos.}
\label{septop}

\blockn{In this subsection, we will first show that for a decidable bound $B$ of an atomic topos $\topos$, the hypergroupoid of atoms of $B \times B$ is locally finite if and only if the slice topos $\topos_{/B}$ is separated (or Hausdorff) in the sense of \cite{moerdijk2000proper}. Then we will show that an atomic topos admits such a bound if and only if it is locally decidable and ``locally separated'', that is if there exists an inhabited object $X$ of $\topos$ such that $\topos_{/X}$ is separated. And finally, that in this case the time evolution constructed in \ref{Qsyst} is completely canonical when seen as a family of functors on the category of Hilbert space of $\topos$  and is described by a canonical principal $\mathbb{Q}_+^*$ bundle $\chi_{\topos} : \topos \rightarrow B\mathbb{Q}_+^*$ attached to every locally separated (locally decidable) atomic topos $\topos$.
}

\block{We recall that a geometric morphism $f: \topos \rightarrow \topos[1]$ is said to be proper if $f_*(\soc)$ is a compact locale internally in $\topos[1]$, and is said to be separated if the diagonal map $\topos \rightarrow \topos \times_{\topos[1]} \topos$ is proper. A topos is said to be compact (resp. separated) if the geometric morphism from $\topos$ to $\set$ is proper (resp. separated). These notions have been studied in \cite{moerdijk2000proper} and in \cite{sketches}C3.2 (see also C5.1).

\bigskip

We will say that a topos is locally separated if there exists an inhabited object $X$ of $\topos$ such that the slice topos $\topos_{/X}$ is separated. }

\block{\label{septop_sepfinit}We start by the following proposition which relates finiteness conditions to the separation property.

\Prop{An atomic locally decidable topos $\topos$ is separated if and only if every atom of $\topos$ is internally finite. }
Also, the ``only if'' part holds without assuming that $\topos$ is locally decidable.

\Dem{We start by assuming that $\topos$ is separated, and that $a \in |\topos|$ is an atom. Then  the topos $\topos_{/a}$ is hyperconnected\footnote{This mean that the locale $p_*(\soc[\topos_{/a}])$ whose open are subobjects of $a$ is trivial, which is the definition of the fact that $a$ is an atom.} and hence proper. Proposition II.2.1(iv) of \cite{moerdijk2000proper} asserts that when one has a commutative diagrame:

\[
\begin{tikzcd}[ampersand replacement=\&]
\topos_{/a} \arrow{r}{g} \arrow{rd}{h} \& \topos \arrow{d}{f} \\
 \& * 
\end{tikzcd}
\]
with $h$ proper and $f$ separated then $g$ is proper. But the map $\topos_{/a} \rightarrow \topos$ is proper if and only if the discrete space $a$ (internally in $\topos$) is compact if and only if $a$ is finite.

\bigskip

Conversely let $\topos$ be an atomic topos whose atoms are internally finite.

The commutative diagram:

\[
\begin{tikzcd}[ampersand replacement=\&]
\topos \arrow{r}{\Delta} \arrow{rd}{Id} \& \topos \times \topos \arrow{d}{\pi_2} \\
 \& \topos
\end{tikzcd}
\]

can be seen as a point $\Delta$ of the topos $\topos \times \topos$ internally in $\topos$. As $\topos \times \topos $ is the pullback of $\topos$ by the canonical geometric morphism from $\topos$ to the point, it will still be an atomic topos internally in $\topos$, and it will still have a generating family of finite objects and hence all its atoms will be finite internally in $\topos$.

Hence our problem is equivalent to prove (constructively) that if $\topos$ is atomic with a point $p$ and that all the atoms of $\topos$ are finite then $p$ is proper. But an atomic topos with a point is equivalent to $BG$ for $G$ the localic group of automorphisms of the point, and the fact that the atoms are finite means that all the $G$-transitive sets are finite, and as $G$ has been taken to be set of automorphisms of the point, this implies that $G$ is compact:

Indeed, the localic monoid of endomorphisms $M(G) = \lim G/U$ constructed in \cite{moerdijk1987morita} is compact by (localic) Tychonoff's theorem, and separated because thanks to the locale decidability one can restrict to the $U$ such that $G/U$ is decidable, hence separated. In particular the point $1$ is closed, and as $G$ can be identified with the subspace of $M(G) \times M(G)$ of $f,g$ such that $fg=1$ which is a closed subspace in a separated compact space, $G$ is also compact.

Finally, the map $1 \rightarrow BG$ is proper because its pull-back along itself is the map $G \rightarrow *$ which we just showed to be proper, and the map $* \rightarrow BG$ is always an open surjection (for exemple by \cite{sketches}C3.5.6(i)) hence the fact that proper map descend along open surjection (see \cite{sketches}C5.1.7) allows us to conclude.
}

}

\block{\label{septop_finitdeg}If $e$ and $e'$ are two decidable atoms and $f:e \rightarrow e'$ then the internal (semi-continuous) number $|f^{-1}(y)|$ does not depend on $y$ because $e'$ is an atom and hence gives an (externally) well defined number called the degree of $f$.

\Prop{Let $\topos$ be an atomic topos, and $B$ be a decidable bound of $\topos$ then the associated semi-simple hypergroupoid is locally finite if and only if for any couple $e$,$e'$ of atoms of $B$, any map $f:e \rightarrow e'$ has finite degree.}

\Dem{If $\Gamma_f$ denotes the graph of $f$ in $B \times B$ then it is an atom of $B \times B$ and the degree of $f$ is equal to $|\Gamma_f|_r$. Hence if the hypergroupoid is locally finite, then every such map has a finite degree. Conversely, if any such map has a finite degree then every simple element is right (and left) finite and as any element can be written $g=uv^*$ with $u$ and $v$ simples, one has:

\[ |g|_l = |\{z |\exists t, z=u(t) \text{ and } v(t)=y \}| \]

and this set is finite because it is a quotient of $v^{-1}(y)$.

}

}

\block{\Prop{Let $\topos$ be an atomic topos, $B$ a decidable bound of $\topos$ then the associated hypergroupoid is locally finite if and only if $\topos_{/B}$ is separated. Moreover, for an arbitrary atomic topos, such a bound exists if and only if it is locally decidable and locally separated.}

\Dem{Let $\topos$ be an atomic topos, and $B$ an arbitrary object of $\topos$. One has:

\[ \topos_{/B} = \coprod_{a \text{ atom of } B} \topos_{/a} \]

Hence $\topos_{/B}$ is separated if and only if for every atom $a$ of $B$ the topos $\topos_{/a}$ is separated. Also, by \ref{septop_sepfinit}, $\topos_{/B}$ is separated if and only if for every atom $v$ of $\topos$, every map $v \rightarrow B$ has finite fiber.

In particular, if $B$ is decidable and $\topos_{/B}$ is separated, then for any atom $e$ of $\topos$ and any atom $e'$ of $B$, any map $f:e' \rightarrow e$ has finite degree, and hence the associated hypergroupoid is locally finite.

\smallskip

Conversely, let $B$ be a decidable bound satisfying the condition of \ref{septop_finitdeg}, then in order to show that $\topos_{/B}$ is separated we need to show that for any map $f: e' \rightarrow e$ where $e$ is an atom of $B$ and $e'$ an arbitrary atom of $\topos$, $f$ has finite degree. But in this situation, as $B$ is a bound there exists an atom $e''$ of $B$ and a map $e'' \rightarrow e' \rightarrow e$. As the map $e'' \rightarrow e'$ is surjective because they are atoms, any fiber of the map $e' \rightarrow e$ is covered by a fiber of the map $e'' \rightarrow e$, which are finite by assumption. This concludes the proof of the first part of the theorem.

\bigskip

If $\topos$ admits a bound satisfying the locale finiteness assumption, then it is in particular locally separated. Conversely assume that $\topos$ is atomic, locally decidable and locally separated. Let $X$ be an inhabited object such that $\topos_{/X}$ is separated, $B_0$ an arbitrary bound and $B$ be a decidable cover of $B \times X$. Then $B$ is decidable, it is a bound because it is a cover of $B$ and  $\topos_{/B}$ is separated because $B$ can be seen as a decidable object of $\topos_{/X}$, hence the geometric morphism $\topos_{/B} \rightarrow \topos_{/X}$ is separated (see \cite{moerdijk2000proper} II.1.3(1) ) and by composition of separated morphisms $\topos_{/B}$ is separated (see \cite{moerdijk2000proper}II.2.1(ii) ).
}

}

\block{The time evolution constructed in \ref{Qset} can be seen as a family of functors $\mathcal{H} \rightarrow \mathcal{H}^t$ for $t\in \mathbb{R}$ acting on the category of Hilbert spaces of $\topos$, corresponding to the functor which send a representation of $\Cmax(Q)$ to the representation twisted by $\sigma_t$.

By the previous theorem we know that any atomic locally decidable locally separated topos has such a time evolution. We will show that this time evolution is canonical by giving a construction of it which does not depend on the choice of the bound $B$.

\bigskip

To be more precise, let $\topos$ be a locally decidable locally separated topos. We will construct a $\Qt$ principal bundle $\chi_{\topos}$ in $\topos$ the following way. The decidable atoms $a$ of $\topos$ such that $\topos_{/a}$ is separated form a generating family. Hence to define an object of $\topos$ it is enough to define a sheaf for the atomic topology on the full subcategory of these atoms.

We define : 
\[\hom(a,\chi_{\topos}) = \mathbb{Q}_+^* \]
and if $f:a \rightarrow a'$ is any map then it acts on $\mathbb{Q}_+^*$ by multiplication by its degree. All these maps are bijective, hence it defines a sheaf.  Also $\mathbb{Q}_+^*$ acts on $\chi_{\topos}$ by multiplication turning it into a principal $\mathbb{Q}_+^*$ bundle.

We note in particular that if $\topos$ is itself separated, then the terminal object of $\topos$ is among the atoms of the site we consider and hence $\hom(1,\chi_{\topos})$ is inhabited and hence $\chi_{\topos}$ is the trivial bundle. But saying that $\chi_{\topos}$ is trivial only means that it is possible to construct a global section $d$ of it, which will be a map associating to any decidable atom $a$ such that $\topos_{/a}$ is separated a finite rational number $d(a)$ such that if $f:a \rightarrow a'$ is a map between two of these atoms then $d(a)=d(a')deg(f)$. 

}

\block{Finally, the time evolution given by any bound can be described in terms of this invariant $\chi_{\topos}$. Indeed if $\mathcal{H}$ is an arbitrary Hilbert space on $\topos$, and  we choose an 'admissible' bound $B$, then the effect of the time evolution on $\mathcal{H}$ can be described by the fact that $\mathcal{H}^t$ is the same sheaf as $\mathcal{H}$ on the site of atoms of $B$ but with the action of a map $f:e\rightarrow e'$ twisted by $\chi(f)^{it} = (deg f) ^ {-it}$.

Hence if we see $\chi_{\topos}$ as a morphism from the topos $\topos$ to the classifying space$B\mathbb{Q}_+^*$ of principal $\mathbb{Q}_+^*$ bundles (i.e. the topos of $\mathbb{Q}_+^*$ sets), and if we call $E_t$ the one dimensional Hilbert space in $B\mathbb{Q}_+^*$ defined by $\mathbb{C}$ with its usual Hilbert space structure and endowed with action $q.z = q^{-it}(z)$ then the previous formula for $\mathcal{H}^t$ can be rephrased as:

\[ \mathcal{H}^t = \mathcal{H} \otimes \chi_{\topos}^*(E_t) \]

where the tensor product is just the internal tensor product of Hilbert spaces in $\topos$.
}

\subsection{Examples}

\block{If $G$ is a discrete group, $\topos$ the topos of $G$-sets and $B$ is $G$ endowed with its (left) action on itself, then the corresponding quantale is $\mathcal{P}(G)$, the hypergroupoid is the group $G$, the integral algebra is $\mathbb{Z}[G]$, and the reduced and maximal $C^*$-algebras are the usual reduced and maximal group $C^*$ algebras. In this situation the time evolution is trivial.}

\block{The best-known example of this situation is the case of double cosets algebras. Let $G$ be a discrete group, and $(K_i)$ a family of subgroups of $G$ (one can generalize to $G$ a localic group and $K_i$ open subgroups). Let $X_i$ be the $G$-set $G/K_i$. The topos of $G$-set is atomic and $\rel(X_i,X_j)$ can be identified with the subset of $G$ stable by the action of $K_i$ on the left and $K_j$ one the right, hence the atom of $\rel(X_i,X_j)$ are exactly the $(K_i,K_j)$ double-cosets. Under this identification, the composition of a $(K_i,K_j)$ cosets with a $(K_j,K_t)$ cosets is the set of $(K_i,K_t)$ cosets included in the product of their elements, and the coefficients $\hgacoef{a}{g}{g'}$ are exactly the usual coefficient involved in the definition of the double cosets modules and double cosets algebras (hence we are in the locally finite case if and only if one has the usual almost normality condition). }

\block{The previous example in particular gives back the BC-system constructed in \cite{bost1995hecke} with both its time evolution and its integral sub-algebra by considering the topos of continuous actions of the group $G = \mathbb{A}^f_\mathbb{Q} \rtimes \mathbb{Q^*_+}$, where $\mathbb{A}^f_\mathbb{Q}$ denote the additive group of finite adele of $\mathbb{Q}$, i.e. the restricted product of all $p$-adic completions of $\mathbb{Q}$, with the bound $B=G/\hat{\mathbb{Z}}$.

Unfortunately, trying to replace $\mathbb{Q}$ by another number field in this construction does not  seem to give the ``good'' $BC$-system for number field constructed in \cite{ha2005bost}, and certainly not the good arithmetic subalgebra. Actually the variant of the BC algebra associated to a number field $K$ constructed in \cite{arledge1997semigroup} corresponds to the topos of continuous $G_K=(\mathbb{A}_K^f) \ltimes K^*$-sets with the bound $B = G_K / \widehat{\mathcal{O}_K}$. }

\block{As a spatial atomic topos is automatically a disjoint sum of toposes of the form $BG$ for $G$ a localic group, the case of double cosets algebras (attached to discrete, topological or localic groups with possibly several subgroups involved) exactly corresponds to the case of spatial atomic toposes. We actually do not know if in the classical set theoretical case there exists examples of non-spatial locally separated atomic toposes, but these examples definitely exist internally to other toposes, and we will exploit them in a forthcoming article.   }

\bibliography{Biblio}{}
\bibliographystyle{plain}

\renewcommand{\thefootnote}{\fnsymbol{footnote}} 
\setcounter{footnote}{0}
\footnotetext{\textsc{Simon Henry, Paris, France.}}
\footnotetext{\emph{E-mail adress:} \texttt{henrys@math.jussieu.fr}}   
\renewcommand{\thefootnote}{\arabic{footnote}} 

\end{document}